\documentclass[12pt]{amsart}
\usepackage{amsmath}
\usepackage{amssymb}
\usepackage{latexsym}
\usepackage{amscd}
\usepackage{pgf,tikz}

\newdimen\AAdi%
\newbox\AAbo%
\def\AAk#1#2{\s_etbox\AAbo=\hbox{#2}\AAdi=\wd\AAbo\kern#1\AAdi{}}%
\def\AAr#1#2#3{\s_etbox\AAbo=\hbox{#2}\AAdi=\ht\AAbo\raise#1\AAdi\hbox{#3}}%
\font\tenmsb=msbm10 at 12pt
\font\sevenmsb=msbm7 at 8pt
\font\fivemsb=msbm5 at 6pt
\newfam\msbfam
\textfont\msbfam=\tenmsb
\scriptfont\msbfam=\sevenmsb
\scriptscriptfont\msbfam=\fivemsb

\newcommand{\ba}{\begin{array}}
\newcommand{\ea}{\end{array}}

\numberwithin{equation}{section}

\begin{document}
\title
[Inequalities on Finsler Manifolds] {Generalized Hardy Type and Caffarelli-Kohn-Nirenberg Type Inequalities on Finsler Manifolds}
\author
[Shihshu Walter Wei and Bingye Wu]{Shihshu Walter Wei$^*$ and Bing Ye Wu$^{**}$}
\address[Shihshu Walter Wei] 
{Department of Mathematics, University of Oklahoma, Norman, Oklahoma, 73019-0315, U.S.A. }
\email{wwei@ou.edu}

\address[Bing-Ye Wu]
{Department of Mathematics, Minjiang University,
Fuzhou, Fujian, 350108, China} \email{wubingye@mju.edu.cn}

\date{}
\maketitle

\begin{abstract}
In this paper we derive both local and global geometric inequalities on general Riemannnian and Finsler manifolds and prove generalized Caffarelli-Kohn-Nirenberg type and Hardy type inequalities on Finsler manifolds, illuminating curvatures of both Riemannian and Finsler manifolds influence geometric inequalities.
\end{abstract}



\section{Introduction}

It is well-known that Hardy type inequalities have been widely used in analysis and differential equations. In \cite{CKN} Caffarelli, Kohn and Nirenberg proved rather general interpolation inequalities with weights. Recently in \cite{WL}, Wei and Li used comparison theorems in Riemannian geometry to prove some sharp generalized Hardy type and Caffarelli-Kohn-Nirenberg type inequalities on Riemannian manifolds.
Some applications of generalized Hardy type inequalities in $p$-harmonic geometry have been studied in \cite{CLW1}.

Finsler geometry, as the natural generalization of Riemannian geometry, has been a very active field in differential geometry and appears in a broad spectrum of contexts. (e.g., two different Finsler metrics, the Kobayashi metric and the Carath\'eodory metric appear very naturally in the theory of several complex variables.) The main purpose of the present paper is to, on the one hand, give a local and two $L^p$ versions of the results in \cite{WL} so that the inequalities work on every Riemannian manifold and in a wider class, and on the other hand,  generalize their results from Riemannian manifolds to Finsler manifolds. We use Hessian and Laplacian comparison theorems in Finsler Geometry by constructing appropriate vector fields. It should be pointed out here that
the volume form on a Riemannian manifold is uniquely determined by the given Riemannian metric,
while there are different choices of volume forms for Finsler
metrics. The frequently used volume forms in Finsler geometry are the so-called Busemann-Hausdorff volume form and Holmes-Thompson volume form, and in \cite{Wu1,Wu2} we introduce the extreme volume forms $dV_{\rm ext}$ (include the maximal and minimal volume forms, cf. \eqref{2.8} and \eqref{2.9}) for Finsler manifolds
which also play the important role in Finsler geometry. In this paper we shall mainly use the extreme volume forms.

To state our results we need some notions from Finsler geometry, for details see \S 2. Throughout this paper, unless otherwise stated, we let $(M,F)$ be a complete Finsler manifold with finite uniformity constant $\mu_F$ (cf. \eqref{2.11}),  Cut($x_0$) be the cut locus of a fixed point $x_0$, and $\Omega\subset M\backslash{\rm Cut}(x_0)$ be a domain in $M.$ It should be pointed out here that in general there are three completeness for Finsler manifolds: forward complete, backward complete and complete (i.e., both forward and backward complete), and they are equivalent when $\mu_F<\infty$. In this situation, the distance function $r=d_F(x_0,\cdot):\Omega\rightarrow\mathbb{R}$  from $x_0$ is smooth on $\Omega\backslash\{x_0\}$, and thus the gradient vector field $\nabla r$ of $r$ (with respect to Finsler metric $F$) is also smooth on $\Omega\backslash\{x_0\}$. We usually call  $\nabla r$ the {\it radial vector field with respect to} $x_0$. We call $x_0\in M$ a {\it pole}, if the exponential map exp$_{x_0}:T_xM\rightarrow M$ is a diffeomorphism. We say that $M$ {\it has nonpositive} (resp. {\it nonnegative}) {\it radial flag curvature at} $x_0$ if flag curvature {\bf K}$(\nabla r;P)$ of flag $(\nabla r;P)$ whose flag pole is a radial vector is nonpositive (resp. nonnegative) for every plane $P$ (cf. \eqref{2.5}). Similarly, we say that $M$ {\it has nonpositive} (resp. {\it nonnegative}){\it radial Ricci curvature} {\bf Ric} $(\nabla r) $ at $x_0$ if  {\bf Ric} $(\nabla r) \leqslant0$ (resp. $\geqslant0$) (cf. \eqref{2.6}).
In this paper we first derive both local and global
{\bf Geometric Inequalities 4.1 and 4.2} on every Riemannnian manifold and Finsler manifold $\big ($cf. $(4.1.a), (4.1.b), (4.2.a),  (4.2.b)\big ).$
We then prove generalized Caffarelli-Kohn-Nirenberg type and Hardy type inequalities on Finsler manifolds.
The main results of this paper are the following:

\noindent
{\bf Theorem 1.1} {\it Let $(M,F)$ be an $n$-dimensional complete Finsler manifold with finite uniformity constant $\mu_F$. Let ${\rm Cut}(x_0)$ be the cut locus of a fixed point $x_0$, and $\Omega\subset M\backslash{\rm Cut}(x_0)$ be a domain in $M\, .$ Suppose that the radial flag curvature ${\bf K}(\nabla r;\,\cdot)$
or radial Ricci curvature ${\bf Ric}(\nabla r)$ of $M$ satisfies one of the following three conditions:

$($i$)$ $0\leqslant {\bf Ric}(\nabla r)$ and $n \leqslant a+b+1$;

$($ii$)$ ${\bf K}(\nabla r;\,\cdot) \leqslant 0$ and $a+b+1\leqslant n$;

$($iii$)$  ${\bf K}(\nabla r;\,\cdot)  = 0$ and $a,b\in
\mathbb{R}$ are any constants.

Then for any $u\in C_0^{\infty}(\Omega\backslash \{x_0\})$, the following Caffarelli-Kohn-Nirenberg type inequality holds:
\begin{equation}\label{11}
\int_{\Omega}\frac{|u|^2}{r^{a+b+1}}dV_{\rm ext}\leqslant\hat{\mu}_F^{\frac{n+1}{2}}\cdot
\left(\int_{\Omega}\frac{|u|^{p}}{r^{ap}}dV_{\rm ext}\right)^{\frac {1}{p}}\left(\int_{\Omega}\frac{(F(\nabla u))^{q}}{r^{bq}}dV_{\rm ext}\right)^{\frac {1}{q}}.
\end{equation} In particular, if $M$ has a pole $x_0$ or $\operatorname{Cut}(x_0)$ is empty in $($i$)$, or $M$ is simply connected in $($ii$)$ or $($iii$)$, then for any $u\in C_0^{\infty}(M\backslash \{x_0\})$,
\begin{equation}\label{12}
\int_{M}\frac{|u|^2}{r^{a+b+1}}dV_{\rm ext}\leqslant\hat{\mu}_F^{\frac{n+1}{2}}\cdot
\left(\int_{M}\frac{|u|^{p}}{r^{ap}}dV_{\rm ext}\right)^{\frac {1}{p}}\left(\int_{M}\frac{(F(\nabla u))^{q}}{r^{bq}}dV_{\rm ext}\right)^{\frac {1}{q}},
\end{equation}
where $\hat{\mu}_F^{\frac{n+1}{2}}= \mu_F^{\frac{n+1}{2}} \cdot |\frac{2}{n-a-b-1}\big |$,  $p\in [1, \infty]$, $\frac 1p + \frac 1q = 1,$ and if $p=\infty$, $\left(\int_{M}\frac{|u|^{p}}{r^{ap}}dV_{\rm ext}\right)^{\frac {1}{p}}$ stands for the supreme of $\frac{|u|}{r^{a}}.$}

This result is new, even when $M$ is a Riemannian manifold:

\noindent
{\bf Corollary 1.1} {\it Let $M$ be a complete $n$-dimensional Riemanian manifold with the volume element $dv$. Fix $x_0\in M$ let $\Omega\subset M\backslash{\rm Cut}(x_0)$ be a domain in $M\, .$ Suppose that radial curvature $K_r$, or radial Ricci curvature $\text {Ric}_{rad}$ of $\Omega$ satisfies  one of the following three conditions:

$($i$)$ $0\leqslant \text {Ric}_{rad}$ and $n \leqslant a+b+1$;

$($ii$)$ $K_r \leqslant 0$ and $a+b+1\leqslant n$;

$($iii$)$  $K_r  = 0$ and $a,b\in
\mathbb{R}$ are any constants.

Then for any $u\in C_0^{\infty}(\Omega\backslash\{x_0\})$, the following Caffarelli-Kohn-Nirenberg type inequality holds:
\begin{equation}\label{13}
C\cdot\int_{\Omega}\frac{|u|^2}{r^{a+b+1}}dv \leqslant
\left(\int_{\Omega}\frac{|u|^{p}}{r^{ap}}dv\right)^{\frac {1}{p}}\left(\int_{\Omega}\frac {|\nabla u|^{q}}{r^{bq}}dv\right)^{\frac {1}{q}}.
\end{equation} In particular, if in $($i$)$ $\operatorname{Cut}(x_0) = \emptyset$, or in $($ii$)$ or in $($iii$)$ $\pi_1(M)=0$, then for any $u\in C_0^{\infty}(M\backslash \{x_0\})$,
\begin{equation}\label{14}
C\cdot\int_{M}\frac{|u|^2}{r^{a+b+1}}dv \leqslant
\left(\int_{M}\frac{|u|^{p}}{r^{ap}}dv\right)^{\frac {1}{p}}\left(\int_{M}\frac {|\nabla u|^{q}}{r^{bq}}dv\right)^{\frac {1}{q}},
\end{equation} where $C=C(a,b)=\big |\frac{n-a-b-1}{2}\big |$,  $p\in [1, \infty]$, and $\frac 1p + \frac 1q = 1.$}

\noindent
{\bf Theorem 1.2} {\it Let $(M,F), \Omega$ and curvature conditions as in Theorem 1.1 $(i), (ii), (iii)$ Then for any $u\in C_0^\infty(\Omega\backslash\{x_0\})$, the following Caffarelli-Kohn-Nirenberg type inequality holds:
\begin{equation}\label{15}
\int_{\Omega}\frac{|u|^p}{r^{a+b+1}}dV_{\rm ext}\leqslant\check{\mu}_F^{\frac{n+1}{2}}\cdot
\left(\int_{\Omega}\frac{|u|^{p}}{r^{aq}}dV_{\rm ext}\right)^{\frac {1}{q}}\left(\int_{\Omega}\frac{(F(\nabla u))^{p}}{r^{bp}}dV_{\rm ext}\right)^{\frac {1}{p}}.
\end{equation} In particular, if $M$ has a pole $x_0$ or $\operatorname{Cut}(x_0)$ is empty in $($i$)$, or $M$ is simply connected in $($ii$)$ or $($iii$)$, then for any $u\in C_0^\infty(M\backslash\{x_0\})$,
\begin{equation}\label{16}
\int_{M}\frac{|u|^p}{r^{a+b+1}}dV_{\rm ext}\leqslant\check{\mu}_F^{\frac{n+1}{2}}\cdot
\left(\int_{M}\frac{|u|^{p}}{r^{aq}}dV_{\rm ext}\right)^{\frac {1}{q}}\left(\int_{M}\frac{(F(\nabla u))^{p}}{r^{bp}}dV_{\rm ext}\right)^{\frac {1}{p}},
\end{equation}
where $\check{\mu}_F^{\frac{n+1}{2}}=\mu_F^{\frac{n+1}{2}}\cdot \big |\frac{p}{n-a-b-1}\big |$,  $p\in (1, \infty)$ and $\frac 1p + \frac 1q = 1.$}

Similarly, when $M$ is a Riemannian manifold, we have the following new result:

\noindent
{\bf Corollary 1.2} {\it Let $M, \Omega$ and curvature conditions as in Corollary 1.1 $(i), (ii), (iii)$ Then for any $u\in C_0^\infty(\Omega\backslash\{x_0\})$, the following Caffarelli-Kohn-Nirenberg type inequality holds:
\begin{equation}\label{17}
\tilde{C}\cdot\int_{\Omega}\frac{|u|^p}{r^{a+b+1}}dv \leqslant
\left(\int_{\Omega}\frac{|u|^{p}}{r^{aq}}dv\right)^{\frac {1}{q}}\left(\int_{\Omega}\frac {|\nabla u|^{p}}{r^{bp}}dv\right)^{\frac {1}{p}}.
\end{equation} In particular, if in $($i$)$ $\operatorname{Cut}(x_0) = \emptyset$, or in $($ii$)$ or in $($iii$)$ $\pi_1(M)=0$, then for any $u\in C_0^\infty(M\backslash\{x_0\})$,
\begin{equation}\label{18}
\tilde{C}\cdot\int_{M}\frac{|u|^p}{r^{a+b+1}}dv \leqslant
\left(\int_{M}\frac{|u|^{p}}{r^{aq}}dv\right)^{\frac {1}{q}}\left(\int_{M}\frac {|\nabla u|^{p}}{r^{bp}}dv\right)^{\frac {1}{p}},
\end{equation} where $\tilde{C}=\tilde{C}(a,b)=\big |\frac{n-a-b-1}{p}\big |$,  $p\in (1, \infty)$ and $\frac 1p + \frac 1q = 1.$}

When $p=q=2,$ the above two Theorems $1.1$ and $1.2$ meet and give their $L^2$ version:

\noindent
{\bf Theorem 1.3} {\it Let $(M,F)$ be an $n$-dimensional  complete Finsler manifold with finite uniformity constant $\mu_F$. Fix $x_0\in M$ let $\Omega\subset M\backslash{\rm Cut}(x_0)$ be a domain. Suppose that the radial flag curvature ${\bf K}(\nabla r;\,\cdot)$
or radial Ricci curvature ${\bf Ric}(\nabla r)$ satisfies one of the following three conditions:

$($i$)$ $0\leqslant {\bf Ric}(\nabla r)$ and $n \leqslant a+b+1$;

$($ii$)$ ${\bf K}(\nabla r;\,\cdot) \leqslant 0$ and $a+b+1\leqslant n$;

$($iii$)$  ${\bf K}(\nabla r;\,\cdot)  = 0$ and $a,b\in
\mathbb{R}$ are any constants.

Then for any $u\in C_0^{\infty}(\Omega\backslash\{x_0\})$, the following Caffarelli-Kohn-Nirenberg type inequality holds:
\begin{equation}\label{19}
\int_{\Omega}\frac{|u|^2}{r^{a+b+1}}dV_{\rm ext}\leqslant\hat{\mu}_F^{\frac{n+1}{2}}\cdot
\left(\int_{\Omega}\frac{|u|^{2}}{r^{2a}}dV_{\rm ext}\right)^{\frac {1}{2}}\left(\int_{\Omega}\frac{(F(\nabla u))^{2}}{r^{2b}}dV_{\rm ext}\right)^{\frac {1}{2}}.
\end{equation} In particular, if $M$ has a pole $x_0$ or $\operatorname{Cut}(x_0)$ is empty in $($i$)$, or $M$ is simply connected in $($ii$)$ or $($iii$)$, then for any $u\in C_0^{\infty}(M\backslash\{x_0\})$,
\begin{equation}\label{110}
\int_{M}\frac{|u|^2}{r^{a+b+1}}dV_{\rm ext}\leqslant\hat{\mu}_F^{\frac{n+1}{2}}\cdot
\left(\int_{M}\frac{|u|^{2}}{r^{2a}}dV_{\rm ext}\right)^{\frac {1}{2}}\left(\int_{M}\frac{(F(\nabla u))^{2}}{r^{2b}}dV_{\rm ext}\right)^{\frac {1}{2}}.
\end{equation}
where $\hat{\mu}_F^{\frac{n+1}{2}}=\mu_F^{\frac{n+1}{2}} \cdot \big |\frac{2}{n-a-b-1}\big |$.}

\noindent
{\bf Remark} When $F$ is Riemannian, one has $\mu_F=1$, and we recapture the corresponding results for Riemannian manifolds \cite{WL, CLW1, CLW2, W}.

As applications, we obtain embedding theorems for weighted Sobolev spaces of functions on Finsler manifolds (cf. Theorem 7.1) and geometric differential-integral inequalities on Finsler manifolds (cf. Theorem 7.2), generalizing the work in \cite{WL} in Riemannian manifolds. We then focus our study on generalized Hardy type inequalities on Finsler manifolds by using the {\it double} limiting technique in \cite {WL}, and extend the density argument in \cite{CLW2}. We introduced the notion of the  space $W_{F,0}^{1,p}(M)$ on a Finsler manifold to be the completion of smooth compactly supported functions $u\in C_0^{\infty}(M)$
with respect to the ``norm"
\begin{equation}\label{1.11}
\|u\|_{W_{F,0}^{1,p}(M)}:= \left(\int_M \big (  |u|^p +
(F(\nabla u))^p  \big )\,  dV_{\rm ext}  \right)^{\frac{1}{p}}.
\end{equation}
It is easy to verify that $\|\cdot\|_{W_{F,0}^{1,p}(M)}$ satisfies the following properties:\\\indent
(i) (Positive definiteness) $\|u\|_{W_{F,0}^{1,p}(M)}\geqslant0, \forall u\in W_{F,0}^{1,p}(M)$, and $\|u\|_{W_{F,0}^{1,p}(M)}=0$ if and only if $u=0$ almost everywhere.\\\indent
(ii) (Positive homogeneity) $\|\lambda u\|_{W_{F,0}^{1,p}(M)}=\lambda\|u\|_{W_{F,0}^{1,p}(M)},\forall \lambda>0$ and $u\in W_{F,0}^{1,p}(M)$.\\\indent
(iii) (Triangle inequality) $\|u+v\|_{W_{F,0}^{1,p}(M)}\leqslant \|u\|_{W_{F,0}^{1,p}(M)}+\|v\|_{W_{F,0}^{1,p}(M)}, \forall u,v\in W_{F,0}^{1,p}(M)$.\\\indent
We note here that $\|\lambda u\|_{W_{F,0}^{1,p}(M)}=|\lambda|\cdot\|u\|_{W_{F,0}^{1,p}(M)}$ does not hold for general Finsler metric, that is to say, $\|\cdot\|_{W_{F,0}^{1,p}(M)}$  is not a genuine norm for general Finsler manifold. Nevertheless we have $\|\lambda u\|_{W_{F,0}^{1,p}(M)}\leqslant\mu_F^\frac12|\lambda|\cdot\|u\|_{W_{F,0}^{1,p}(M)}$, and since we assume $\mu_F<\infty$, we may define Cauchy sequence $\{u_i\}\subset C_0^\infty(M)$ with respect to $\|\cdot\|_{W_{F,0}^{1,p}(M)}$ in the usual way, and thus $W_{F,0}^{1,p}(M)$ is well-defined.
We say {\it $\frac {u}{r} \in L^p(M)$ in a Finsler sense}, denoted by $\frac {u}{r} \in L_{F}^p(M)$ if $\int_M\frac{|u|^p}{r^p}dV_{\rm ext} < \infty\, .$
In particular, we have

\noindent
{\bf Theorem 1.4} {\it Let $(M,F)$ be an $n$-dimensional  complete Finsler manifold with nonpositive radial flag curvature at the pole $x_0\in M$ and with finite uniformity constant $\mu_F$. Then for any $u\in W_{F,0}^{1,p}(M)$ and $1<p<n$, the following Hardy type inequality holds:
\begin{equation}\label{111}
\int_M\frac{|u|^p}{r^p}dV_{\rm ext}\leqslant\acute{\mu}_F^{\frac{n+p}{2}}\cdot \int_M(F(\nabla u))^p \, dV_{\rm ext},
\end{equation}
where $\acute{\mu}_F^{\frac{n+p}{2}} = \mu_F^{\frac{n+p}{2}} \cdot \left(\frac{p}{n-p}\right)^p\, .$ Furthermore, $\frac {u}{r} \in L^p(M)\, $ in a Finsler sense.}

This recaptures a result of Wei-Li \cite[Theorem 1, Corollary 1.2]{WL} (cf. Corollary 8.2), when $M$ is a Riemannian manifold.

\noindent
{\bf Theorem 1.5} {\it Let $(M,F)$ be an $n$-dimensional  complete Finsler manifold with nonnegative radial Ricci curvature at the pole $x_0\in M$, and with finite uniformity constant $\mu_F$. Then for any $u\in W_{F,0}^{1,p}(M),$ $\frac {u}{r} \in L_{F}^p(M)$  and $p>n$, the following Hardy type inequality holds:}
\begin{equation}\label{112}
\int_M\frac{|u|^p}{r^p}dV_{\rm ext}\leqslant\grave{\mu}_F^{\frac{n+p}{2}}\cdot \int_M(F(\nabla u))^p\, dV_{\rm ext},
\end{equation}
where $\grave{\mu}_F^{\frac{n+p}{2}} = \mu_F^{\frac{n+p}{2}} \cdot \left(-\frac{p}{n-p}\right)^p $

This recaptures a theorem of Chen-Li-Wei \cite[Theorem 5]{CLW2} (cf. Corollary 8.4), when $M$ is a Riemannian manifold. Furthermore, the assumption $\frac {u}{r} \in L_{F}^p(M)$ cannot be dropped, or a counter-example is constructed in Section 5 in Chen-Li-Wei \cite{CLW2}.

\noindent
{\bf Corollary 1.3} {\it Let $(M,F)$ be an $n$-dimensional  complete Finsler manifold with vanishing flag curvature at the pole $x_0\in M$, and with finite uniformity constant $\mu_F$. Then $(i)$ for any $u\in W_{F,0}^{1,p}(M),$ $\frac {u}{r} \in L_{F}^p(M)$ and $1 < p < \infty$, or (ii)  for any $u\in W_{F,0}^{1,p}(M),$ and $p < n$ the following Hardy type inequality holds:
\begin{equation}\label{113}
\int_M\frac{|u|^p}{r^p}dV_{\rm ext}\leqslant\dot{\mu}_F^{\frac{n+p}{2}}\cdot \int_M(F(\nabla u))^pdV_{\rm ext},
\end{equation}
where $\dot{\mu}_F^{\frac{n+p}{2}} = \mu_F^{\frac{n+p}{2}} \cdot \left|\frac{p}{n-p}\right|^p\, .$}

Theorem 1.5 is in contrast to Theorem 1.4, in which $\frac {u}{r} \in L^p(M)$ in a Finsler sense is a conclusion, rather than an assumption.
The above theorems  illuminate that curvatures of both Riemannian and Finsler manifolds influence geometric inequalities such as generalized Hardy Type and Caffarelli-Kohn-Nirenberg Type inequalities.

\section{Finsler Geometry}
In this section we shall recall some basic notations and formulas in Finsler geometry, for details we refer to \cite{Wu1,Wu3,WX}.
Let $(M,F)$ be a
Finsler $n$-manifold with Finsler metric $F:TM\rightarrow
[0,\infty)$, where $TM$ is the tangent bundle of $M\, .$ Let $(x,y)=(x^i,y^j)$ be  local coordinates on $TM$,
and $\pi:TM\backslash\{0\}\rightarrow M$ be the natural projection. Unlike
in the Riemannian case, in general Finsler quantities are functions defined on
$TM$ rather than $M$. The {\it fundamental tensor} $g_{ij}$ and the
{\it Cartan tensor} $C_{ijk}\, , 1 \leqslant i, j, k \leqslant n$ are defined by
$$g_{ij}(x,y):=\frac{1}{2}\frac{\partial^2F^2(x,y)}{\partial y^i\partial
y^j},\quad
C_{ijk}(x,y):=\frac{1}{4}\frac{\partial^3F^2(x,y)}{\partial
y^i\partial y^j\partial y^k}.$$
According to \cite{Chern}, the pulled-back bundle $\pi^*TM$ admits a
unique affine connection, called the {\it Chern connection}. Its
connection forms $\omega^i_j$ are characterized by the following two structure equations:

\noindent
$\bullet$ Torsion freeness:
$$dx^j\wedge\omega^i_j=0;$$
$\bullet$ Almost $g$-compatibility:
$$dg_{ij}-g_{kj}\omega^k_i-g_{ik}\omega^k_j=2C_{ijk}(dy^k+N^k_ldx^l),$$
where $N^k_l$ are real-valued functions determined by $N^k_ldx^l=y^l\omega^k_l$.
It is easy to know that torsion freeness is equivalent to the
absence of $dy^k$ terms in $\omega^i_j$; namely,
$$\omega^i_j=\Gamma^i_{jk}dx^k,$$
together with the symmetry
$$\Gamma^i_{jk}=\Gamma^i_{kj}.$$
Let $V=V^i\partial/\partial x^i$ be a non-vanishing vector field on
an open subset $\mathcal{U}\subset$ $M$. One can introduce a
Riemannian metric $\widetilde{g}^{V} (\cdot,\cdot) =\langle\cdot,\cdot\rangle_V$ in the direction of $V\, ,$ and a linear connection
$\nabla^V$ on the tangent bundle over $\mathcal{U}$ as follows:
$$\widetilde{g}^{V} (X,Y)=\langle X,Y\rangle_V:=X^iY^jg_{ij}(x,V),\quad \forall X=X^i\frac{\partial}{\partial
x^i},Y=Y^j\frac{\partial}{\partial x^j};$$
$$\nabla^V_{\frac{\partial}{\partial x^i}}\frac{\partial}{\partial
x^j}:=\Gamma^k_{ij}(x,V)\frac{\partial}{\partial x^k}.$$ From the
torsion freeness and almost $g$-compatibility of Chern connection we
have
\begin{equation}\label{21}
\nabla^V_XY-\nabla^V_YX=[X,Y],
\end{equation}
\begin{equation}\label{22}X\cdot \langle Y,Z\rangle_V=\langle \nabla^V_XY,Z\rangle_V+\langle Y,\nabla^V_XZ\rangle_V+2{\bf C}_V(\nabla^V_XV,Y,Z),
\end{equation}
here ${\bf C}_V$ is defined by
$${\bf C}_V(X,Y,Z)=X^iY^jZ^kC_{ijk}(x,V).$$By definition {\bf C}$_V(X,X,Z)$ is totally symmetric with respect to $X,Y,Z$, and by Euler's lemma  it also satisfies
\begin{equation}\label{23}
{\bf C}_V(V,X,Y)=0.
\end{equation} In view of \eqref{21}-\eqref{23} we see that
the Chern connection $\nabla^V$ and the Levi-Civita connection
$\widetilde{\nabla}^V$ of $\widetilde{g}^{V}(\cdot,\cdot)=\langle\cdot,\cdot\rangle_V$ are related by
\begin{equation}\label{24}
\begin{aligned}
\langle \nabla^V_XY,Z\rangle_V& =\langle\widetilde{\nabla}^V_XY,Z\rangle_V-{\bf
C}_V(\nabla^V_XV,Y,Z)\\
& \quad -{\bf C}_V(\nabla^V_YV,X,Z)+{\bf
C}_V(\nabla^V_ZV,X,Y).
\end{aligned}
\end{equation}
The  {\it Chern curvature}
${\bf R}^V(X,Y)Z$ for vector fields $X,Y,Z$ on $\mathcal{U}$ is
defined by
$${\bf R}^V(X,Y)Z:=\nabla^V_X\nabla^V_YZ-\nabla^V_Y\nabla^V_XZ-\nabla^V_{[X,Y]}Z.$$
In the Riemannian case this curvature does not depend on $V$ and
coincides with the Riemannian curvature tensor. Let $P\subset T_xM$ be a 2-plane and $V\in P$ be a nonzero vector. We call the pair $(V;P)$ a {\it flag with pole} $V$.
The {\it flag
curvature} ${\bf K}(V;P)$ of given flag is defined as follows:
\begin{equation}\label{2.5}{\bf K}(V;P)={\bf K}(V;W):=\frac{\langle{\bf R}^V(V,W)W,V\rangle_V}{\langle V,V\rangle_V\langle W,W\rangle_V-\langle V,W\rangle_V^2} .\end{equation}
Here $W$ is a tangent vector such that $V,W$ span the 2-plane $P$
and $V\in T_xM$ is extended to a geodesic field, i.e.,
$\nabla^V_VV=0$ near $x$. The {\it Ricci curvature}\index{Ricci curvature} ${\bf Ric}(V)$ of $V\in T_xM$ is defined by
\begin{equation}\label{2.6}{\bf Ric}(V)=\sum_{i=1}^{n}{\bf K}(V;E_i),\end{equation}
here $E_1,\cdots,E_n$ is the $\widetilde{g}^V$-orthonormal basis for $T_xM$.

Let $f:M\rightarrow \mathbb{R}$ be a smooth function on $M$. The
{\it gradient} $\nabla f$ of $f$ is defined by
$$ df(X)=\langle\nabla f,X\rangle_{\nabla f},\quad \forall X\in \Gamma(TM)$$ whenever $df\ne0$, and $\nabla f=0$ where $df=0$.
Let $\mathcal{U}=$ $\{ x\in M:\nabla f\, |_x\, \ne 0\}$. We define the
{\it Hessian} $Hess(f)$ of $f$ on $\mathcal{U}$ as follows:
$$Hess(f)(X,Y):=X(Y(f))-\nabla^{\nabla f}_XY(f),\quad \forall X,Y\in
\Gamma(T\mathcal{U}).$$ It is known that $Hess(f)$ is symmetric, and it
can be rewritten as (see \cite{WX})
\begin{equation}\label{25}
Hess(f)(X,Y)=\langle \nabla^{\nabla f}_X\nabla
f,Y\rangle_{\nabla f}.
\end{equation} It should be noted  that the notion of Hessian
defined here is different from that in \cite{Sh1}. In that case
$Hess(f)$ is in fact defined by
$$Hess(f)(X,X)=X(X(f))-\nabla^X_XX(f),$$ and there is no definition for
$Hess(f)(X,Y)$ if $X\ne Y$. The advantage of our definition is that
$Hess(f)$ is a symmetric bilinear form and we can treat it by using the
theory of symmetric matrix.

The following Hessian comparison theorem for distance function $r=d_F(x_0,\cdot)$ from $x_0$ first was proved in \cite{WX} with pointwise curvature bounds, and it is easy to see from the proof that the pointwise curvature bounds can be weakened to radial curvature bounds. More precisely, we have the following:

\noindent
{\bf Theorem 2.1 $($Hessian Comparison Theorem under Radial Curvature Assumptions$)$}
 {\it Let $(M,F)$ be an $n$-dimensional  complete Finsler
manifold, and $x_0\in M$.\\\indent
$(1)$ Suppose that $M$ has nonpositive radial flag curvature ${\bf K}(\nabla r;\,\cdot)\leqslant 0$, then for any tangent vector field $X$ on $M$ the following inequality holds whenever $r$ is smooth:
$$Hess(r)(X,X)\geqslant\frac1r\left(\langle X,X\rangle_{\nabla r}-\left\langle\nabla r,X\right\rangle_{\nabla r}^2\right);$$\indent
$(2)$ Suppose that $M$ has nonnegative radial Ricci curvature ${\bf Ric}(\nabla r)\geqslant 0$, then  the following inequality holds whenever $r$ is smooth:
$$\sum_{i=1}^nHess(r)(E_i,E_i)\leqslant\frac{n-1}r,$$ here $E_1,\cdots,E_n$ is the local $\langle\cdot,\cdot\rangle_{\nabla r}$-orthonormal frame on $M$.}

Define
\begin{equation}\label{2.8}
\begin{aligned}dV_{\rm max}&=\sigma_{\rm max}(x)dx^1\wedge\cdots\wedge dx^n\\ {\rm and} \qquad
dV_{\rm min}&=\sigma_{\rm min}(x)dx^1\wedge\cdots\wedge dx^n,\\ \operatorname{where}\qquad
\sigma_{\max} (x):&=\max_{y\in T_xM\setminus
\{0\}}\sqrt{\det\left(\left\langle\frac{\partial}{\partial x^i},\frac{\partial}{\partial x^j}\right\rangle_y\right)},\\
 \sigma_{\rm min}(x):&=\min_{y\in
T_xM\setminus \{0\}}\sqrt{\det\left(\left\langle\frac{\partial}{\partial x^i},\frac{\partial}{\partial x^j}\right\rangle_y\right)}.\end{aligned}
\end{equation}  $\big ($Since $F(x,\lambda y) = \lambda F(x,y)$ for $\forall \lambda > 0\, $ $\Rightarrow$ $\left\langle\frac{\partial}{\partial x^i},\frac{\partial}{\partial x^j}\right\rangle_y := g_{ij}(x,y):=\frac{1}{2}\frac{\partial^2F^2(x,y)}{\partial y^i\partial
y^j}$ is homogeneous of degree $0$ in $y$, both $\max_{y\in T_xM\setminus
\{0\}}$ and $\min_{y\in T_xM\setminus
\{0\}}$ in \eqref{2.8} are taken over a compact set and hence can be realized.$\big)$ Then it is easy to check
that the $n$-forms $dV_{\rm max}$ and $dV_{\rm min}$ as well as the function
$\frac{\sigma_{\rm max}}{\sigma_{\rm min}}$ are well-defined on $M$.
We call $dV_{\rm max}$ and $dV_{\rm min}$ the {\it maximal volume form}
and the {\it minimal volume form} of $(M,F)$, respectively\cite{Wu1,Wu3}. Both
maximal volume form and minimal volume form are called {\it extreme
volume form}, and we shall denote by
\begin{equation}\label{2.9}dV_{\rm ext}\, \, \,  = \quad \rm{either} \quad {\it d}V_{\rm max}\quad \rm{or}\quad {\it d}V_{\rm min}\end{equation}

The {\it uniformity function} $\mu:M\rightarrow\mathbb{R}$ is
defined by
\begin{equation}\label{2.10}
\begin{aligned}\mu(x)&=\max_{y,z,w\in T_xM\backslash\{0\}}\frac{\langle w,w\rangle_y}{\langle w,w\rangle_z}.\end{aligned}\end{equation}
Then
\begin{equation}\label{2.11}
\begin{aligned}\mu_F&=\sup_{x\in M}\mu(x)\end{aligned}\end{equation}
is called the {\it
uniformity constant}.  We always assume $\mu_F<\infty$ throughout this paper. It is known that
\begin{equation}\label{26}
\mu^{-1}_FF^2(w)\leqslant\langle w,w\rangle_y\leqslant\mu_F F^2(w).
\end{equation} Furthermore, we have \cite{Wu1,Wu2}
$$\frac{\sigma_{\rm max}}{\sigma_{\rm min}}\leqslant\mu_F^{\frac{n}{2}},  $$and thus
\begin{equation}\label{27}dV_{\rm max}\leqslant\mu_F^{\frac{n}{2}}dV_{\rm min}.\end{equation}

\section{The Induced Riemannian Metric $\widetilde{g}$ in the Radial Direction $T$}

In this section we introduce the induced Riemannian metric $\widetilde{g}$ by the Finsler metric $F$ in the radial direction $T$. This will play an important role in the proof of main theorems. Let $(M,F)$ be a  complete Finsler manifold, and $x_0\in M$. Then the distance function $r=r(x)=d_F(x_0,x)$ is smooth on $\breve{M}$, where \begin{equation}\breve{M}:=M\backslash\big \lbrace \{x_0\}\cup{\rm Cut}(x_0) \big \rbrace\label{3.0}\end{equation}
and the radial vector field $T=\nabla r$  is also smooth on $\breve{M}$. It is also well-known that $T$ is the unit geodesic field, i.e., $F(T)=1$ and $\nabla^T_TT=0$. Now we define \begin{equation}\label{3.1}\widetilde{g}(\cdot,\cdot) = \langle\cdot,\cdot\rangle_T\, .\end{equation} Then $\widetilde{g}$  is the  Riemannian metric on $\breve{M}$ induced by $F$ in the radial direction $T$. It is clear from the definition of gradient that the gradient $\widetilde{\nabla}r$ of $r$ with respect to $\widetilde{g}$ is just the radial vector field, i.e., $$\widetilde{\nabla}r=\nabla r=T.$$ Thus by \eqref{23},\eqref{24} and \eqref{25} the Laplacian $\widetilde{\Delta}r$ of $r$ with respect to $\widetilde{g}$ is
\begin{equation}
\begin{aligned}{\widetilde{\Delta}}r= & \widetilde{\rm div}(\widetilde{\nabla}r)=\sum_{i=1}^n\widetilde{Hess}(r)(E_i,E_i)=\sum_{i=1}^n\langle\widetilde{\nabla}_{E_i}T,E_i\rangle_T\\
\label{31}= &\sum_{i=1}^n\langle\nabla^T_{E_i}T,E_i\rangle_T=\sum_{i=1}^n Hess(r)(E_i,E_i),
\end{aligned}
\end{equation}
here $\widetilde{\rm div}$ and $\widetilde{Hess}$ are the divergence and Hessian with respect to $\widetilde{g}$, and $E_1,\cdots,E_n$ is the local $\widetilde{g}$-orthonormal field on $\breve{M}$. Therefore, by Theorem 2.1 and \eqref{31} we clearly have the following

\noindent
{\bf Theorem 3.1 (Laplacian Comparison Theorem) } {\it With the same notations as above.\\ \indent
$(1)$ If $M$ has nonpositive radial flag curvature at $x_0$, then $\widetilde{\Delta}r\geqslant\frac{n-1}{r}$.\\ \indent
$(2)$ If $M$ has nonnegative radial Ricci curvature at $x_0$, then $\widetilde{\Delta}r\leqslant\frac{n-1}{r}$.}

Let $u:M\rightarrow\mathbb{R}$ be a smooth function on $M$. By definition, the gradient $\nabla u$ of $u$ (with respect to $F$) and the gradient $\widetilde{\nabla}u$ of $u$ with respect to $\widetilde{g}$ are related by
$$du(X)=\langle\nabla u,X\rangle_{\nabla u}=\langle\widetilde{\nabla} u,X\rangle_{T},\quad\forall X\in\Gamma(T\breve{M}),$$ which together with \eqref{26} and Schwartz inequality yields
$$\langle\widetilde{\nabla}u,\widetilde{\nabla}u\rangle_T=\langle\widetilde{\nabla}u,\nabla u\rangle_{\nabla u}\leqslant F(\nabla u)F(\widetilde{\nabla}u)\leqslant\mu_F^{\frac{1}{2}}F(\nabla u)\langle\widetilde{\nabla}u,\widetilde{\nabla}u\rangle_T^{\frac{1}{2}},$$namely,
\begin{equation}\label{32}|\widetilde{\nabla}u|_{\widetilde{g}}\leqslant\mu_F^{\frac{1}{2}}F(\nabla u).
\end{equation}Here $|\cdot|_{\widetilde{g}}=\langle\cdot,\cdot\rangle_T^{\frac{1}{2}}$ denotes the norm of vector field with respect to $\widetilde{g}$. Let $dV_{\widetilde{g}}$ be the Riemannian volume form of $\widetilde{g}$. It is clear from \eqref{27} that on $\breve{M}$ we have
\begin{equation}\label{33}dV_{\rm min}\leqslant dV_{ \widetilde{g}}\leqslant dV_{\rm max}\leqslant\mu_F^{\frac{n}{2}}dV_{\rm min}.
\end{equation}
\smallskip

A positive {\it Radon measure} $\nu$ is a linear functional on the space $C_0(\mathcal H)$ of real-valued continuous functions on a locally compact Hausdorff space $\mathcal H$ with compact support, $\nu: f \mapsto \nu (f) \in \mathbb R\, ,$ such that $\nu(f) \geqslant 0$, for any $f \geqslant 0$. We define an {\it upper integral} for the nonnegative functions as folloows. 
If $\xi \geqslant 0$ is lower semicontinuous:   $$\nu^{\ast}(\xi) = \sup \nu (f)\, \operatorname{for}\, \operatorname{all}\, \operatorname{nonnegative-valued}\, f \in C_0(\mathcal H)\, \operatorname{satisfying}\, f \leqslant \xi\, ,$$
and for any function $\eta \geqslant 0\, :$ 
$\nu^{\ast}(\eta) = \inf \nu (\xi)\,$ for all lower semicontinuous functions $\xi \, \operatorname{satisfying}\, \eta \leqslant \xi\, .$

A function $f$ is said to be $\nu$-{\it integrable} (or integrable if without ambiguity) if there exists a sequence $\{f_n\} \in C_0(\mathcal H)$ such that $\nu^{\ast}(|f-f_n|) \to 0 $ as $n \to \infty\, .$
A subset $A \subset \mathcal H$ is {\it measurable} with finite measure $\nu (A)\, ,$ if its characteristic function $\chi _A$ is integrable. 
We set $\nu(A) = \int \chi_A\, d\nu\, .$ 

A function $f$ is said to be $\nu$-{\it measurable} (or measurable if without ambiguity) if for all compact sets $K$ and for all $\epsilon > 0\, ,$ there exists a compact set $K_{\epsilon} \subset K\, ,$ such that $\nu (K - K_{\epsilon}) < \epsilon$ and such that the restriction $f \big |_{K_{\epsilon}}$ is continuous on $K_{\epsilon}\, $ (cf. \cite {A}).

\noindent
 {\bf Theorem 3.2.} {\it Let $\breve{M}$ be as in \eqref{3.0}. Then

\noindent
$(1)$ The volume element $dV_{\widetilde{g}}$ on Riemannian manifold $(\breve{M}, \widetilde{g})$  induces a positive Radon measure $\nu$ on a locally compact Hausdorff space $M$.

\noindent
$(2)$ The set $\{x_0\}\cup{\rm Cut}(x_0) $ has measure
$\nu \big ( \{x_0\}\cup{\rm Cut}(x_0) \big )
 = 0$.}

 \begin{proof}
$(1)$ We first note that a Finsler manifold is a locally compact Hausdorff space. For $f \in C_0 (M), $ we define
 \begin{equation}\label{3.5}
 \nu(f) = \int _{M\backslash\lbrace \{x_0\}\cup{\rm Cut}(x_0) \big \rbrace} f\,  dV_{\widetilde{g}} = \int _{\breve{M}} f\,  dV_{\widetilde{g}}.\end{equation}
 Then $\nu$ is a linear functional on $C_0 (M)$ with $\nu(f) \geqslant 0$, for any $f \geqslant 0\, .$ Thus $\nu$ is a positive Radon measure on $M$.

 \noindent
 $(2)$ In view of the definition of the measure of a subset of $M$ and \eqref{3.5}, we have
 the measure $$\begin{aligned}\nu \big ( \{x_0\}\cup{\rm Cut}(x_0) \big ) & = \nu (\chi _{\{x_0\}\cup{\rm Cut}(x_0)}) \\
 & = \int _{\lbrace \{x_0\}\cup{\rm Cut}(x_0) \rbrace \backslash \lbrace \{x_0\}\cup{\rm Cut}(x_0) \rbrace} \chi _{\{x_0\}\cup{\rm Cut}(x_0)} \,  dV_{\widetilde{g}} \\
& = \int _{\emptyset} 0\,  dV_{\widetilde{g}}, \qquad \operatorname{where}\, \emptyset\, \operatorname{is}\, \operatorname{the}\, \operatorname{empty}\, \operatorname{set}\, \\
& = 0 . \end{aligned}$$
 \end{proof}

Assume  $1 \leqslant p \leqslant \infty\, $ for the remaining of this section.

\noindent
{\bf Definition 3.1.}  A measurable function $u : M \to \mathbb R$  is said to {\it belong to $L^p(M)$ with respect to Riemannian metric $\widetilde{g}$ }, denoted by $u \in L^p(M, \widetilde{g})$ if $ \nu(|u|^p)  = \int_{\breve{M}} |u|^p dV_{\widetilde{g}} < \infty\, ,$ and  {\it $u$ is said to belong to $W_0^{1,p}(M)$ with respect to Riemannian metric $\widetilde{g}$}, denoted by $u \in W_0^{1,p}(M, \widetilde{g})$  if there exists a sequence $\{u_i\}$ in $C_0^{\infty}(M)$ such that  $\left(\int_{\breve{M}} |u-u_i|^p + |\widetilde{\nabla} (u- u_i)|^p dV_{\widetilde{g}}\right)^{\frac 1p} \to 0$, as $i \to \infty\, .$
\smallskip

\noindent
 {\bf Theorem 3.3.} {\it $\, L^p(M, \widetilde{g})$ is complete, i.e. every Cauchy sequence $\{u_i\}$ in $L^p(M, \widetilde{g})$ converges $($This means that if for every $\epsilon > 0$, there exists $N$ such that $\big (\nu ^{\ast} (|u_{i} - u_{j}| ^p)\big )^{\frac {1}{p}} = \left(\int_{\breve{M}} |u_{i} - u_{j}| ^p dV_{\widetilde{g}}\right)^{\frac 1p} <  \epsilon$, when $i > N$ and $j > N$, then there exists a unique function $u\in L^p(M, \widetilde{g})$, such that $\big (\nu ^{\ast} (|u_{i} - u| ^p)\big )^{\frac {1}{p}}=\left(\int_{\breve{M}} |u_{i} - u| ^p dV_{\widetilde{g}}\right)^{\frac 1p} \to 0$, as $i \to \infty$ $)$}.

\begin{proof}
In view of Theorem 3.2.(1),  $M$ is a measure space with a Radon measure $\nu$ as defined in \eqref{3.5}. To prove that every Cauchy sequence  in $L^p(M, \widetilde{g})$  converges, it suffices to prove that for every Cauchy sequence $\{u_i\}$ in $L^p(M, \widetilde{g})$, there exists a subsequence $\{u_{i_k}\}$ which converges strongly to a function $u$ in $L^p(M, \widetilde{g})$ as $k \to \infty$, by the triangle inequality. Indeed,
$$\begin{aligned}\left(\int_{\breve{M}} |u_{i} - u| ^p dV_{\widetilde{g}}\right)^{\frac 1p} & \leqslant \left(\int_{\breve{M}} |u_{i} - u_{i_k}| ^p dV_{\widetilde{g}}\right)^{\frac 1p} + \left(\int_{\breve{M}} |u_{i_k} - u| ^p dV_{\widetilde{g}}\right)^{\frac 1p} \\
& < \frac {\epsilon}{2} + \frac {\epsilon}{2}\quad \operatorname{for}\, \operatorname{sufficiently}\, \operatorname{large}\, i, i_k \quad \operatorname{if}\, \{u_{i_k}\}\,  \operatorname{converges}\, \operatorname{to}\, u.
\end{aligned}$$
This also proves that if the limit $u$ exists, then it is unique. The subsequence can be obtained by choosing $u_{i_k}$
so that $\left(\int_{\breve{M}} |u_{i_k} - u_n| ^p dV_{\widetilde{g}}\right)^{\frac 1p} < \frac {1}{2^k}$ for all $n > i_k$. (This is the definition of Cauchy sequence). In particular,  $\left(\int_{\breve{M}} |u_{i_k} - u_{i_{k+1}}| ^p dV_{\widetilde{g}}\right)^{\frac 1p} < \frac {1}{2^k}$ for $k=1,2, \dots$. Furthermore, this subsequence $\{u_{i_k}\}$ gives rise to a bounded monotone sequence of positive functions$$U_{\ell} = |u_{i_1}| +\sum_{k=1}^{\ell}  |u_{i_{k+1}} - u_{i_{k}}| .$$
Indeed, by the triangle inequality $U_{\ell}$ is bounded in $L^p(M, \widetilde{g})$ as
$$\left(\int_{\breve{M}} U_{\ell}^p dV_{\widetilde{g}}\right)^{\frac 1p}  \leqslant \left(\int_{\breve{M}}  |u_{i_1}|^p dV_{\widetilde{g}}\right)^{\frac 1p}  +\sum_{k=1}^{\ell} \frac {1}{2^k}\, = \left(\int_{\breve{M}}  |u_{i_1}|^p dV_{\widetilde{g}}\right)^{\frac 1p}  + \left(1 - \frac {1}{2^{\ell}}\right).$$ By the monotone convergence theorem, $U_{\ell} $ converges pointwise a.e. to a positive function $U$ which is in $L^p(M, \widetilde{g})$ and hence is finite almost everywhere. The sequence $$u_{i_{\ell +1}} = u_{i_1} +\sum_{k=1}^{\ell} \left( u_{i_{k+1}} - u_{i_{k}}\right )$$
thus converges absolutely for almost every $x$ and hence it converges for the same $x$ to some function $u(x)$. Since $|u_{i_k} (x)| \leqslant U_{k-1} (x)\leqslant U(x)$ a.e. and $U \in L^p(M, \widetilde{g})$, by dominated convergence theorem (applying to the sequence $\{|u_{i_k}|^p\}$), $u \in  L^p(M, \widetilde{g})$. Since $ |u_{i_k}(x) - u(x)| \leqslant U(x) + |u(x)| \in
L^p(M, \widetilde{g})\, ,$ again by dominated convergence theorem (applying to the sequence $\{|u_{i_k} - u|^p\}$), we conclude
$$\lim _{k \to \infty}\left(\int _{\breve{M}} |u_{i_k} - u|^p dV_{\widetilde{g}}\right)^{\frac 1p}  = \left(\int _{\breve{M}} (\lim _{k \to \infty}|u_{i_k} - u|)^p dV_{\widetilde{g}}\right)^{\frac 1p} = 0. $$
That is the desired $(\int _{\breve{M}} |u_{i_k} - u|^p dV_{\widetilde{g}})^{\frac 1p} \to 0 $ as $k \to \infty$,  (cf. \cite {LL} for real analysis.)\end{proof}

The proof of Theorem 3.3 yields the following domination and pointwise convergence properties:

\noindent
 {\bf Theorem 3.4.} {\it  If $\{u_i\}$ is a Cauchy sequence in $L^p(M, \widetilde{g})\, ,$ then there exists a subsequence  $\{u_{i_k}\}$ and a nonnegative function $U$ in $L^p(M, \widetilde{g})\, ,$ such that

 \noindent
 $(1)$ $|u_{i_k}| \leqslant U$ almost everywhere in $M .$

 \noindent
 $(2)$ $\lim _{k \to \infty} u_{i_k} = u$ almost everywhere in $M .$}

\section{Geometric Inequalities on General Manifolds}

In this section, we begin with the following geometric inequalities on general Riemannian manifolds and on Finsler Manifolds:

\noindent
{\bf Local and Global Geometric Inequalities 4.1}\label{P:3.1}  {\it \, Let  $\Omega\subset M\backslash{\rm Cut}(x_0)$ be a domain in a Finsler manifold $(M, F).$ For every $u\in C_{0}^{\infty}(\Omega\backslash\{x_0\}),$ and every $a,b\in \mathbb{R},$ the following inequality holds:
\[
 \frac 12 \bigg |\int_{\Omega}\frac{|u|^2}{r^{a+b+1}}\left(r\widetilde{\Delta}r-(a+b)\right)dV_{\widetilde{g}} \bigg |\leqslant \left(\int_{\Omega}\frac{|u|^{p}}{r^{ap}}dV_{\widetilde{g}}\right)^{\frac {1}{p}}
\left(\int_{\Omega}\frac{|\widetilde{\nabla}u|^{q}_{\widetilde{g}}}{r^{bq}}dV_{\widetilde{g}}\right)^{\frac {1}{q}},\tag{4.1.a}
\]
where $p \in [1, \infty]$ with
$\frac 1p+\frac 1q=1$. In particular, if $\operatorname{Cut}(x_0)$
is empty,  or $x_0$ is a pole, then}
\[
 \frac 12 \bigg |\int_{M}\frac{|u|^2}{r^{a+b+1}}\left(r\widetilde{\Delta}r-(a+b)\right)dV_{\widetilde{g}}\bigg |\\
 \leqslant \left(\int_{M}\frac{|u|^{p}}{r^{ap}}dV_{\widetilde{g}}\right)^{\frac {1}{p}}
\left(\int_{M}\frac{|\widetilde{\nabla}u|^{q}_{\widetilde{g}}}{r^{bq}}dV_{\widetilde{g}}\right)^{\frac {1}{q}}.\tag{4.1.b}
\]
\begin{proof} We observe that by (\ref{31}) on any domain $\Omega\subset M\backslash{\rm Cut}(x_0)$,
$$
\begin{aligned}
\widetilde{\rm div}\left(\frac{|u|^2}{r^{a+b}}T\right) & =\frac{|u|^2}{r^{a+b}}\widetilde{\Delta}r+2\widetilde{g}\left(\frac{u}{r^{a+b}}T,\widetilde{\nabla} u\right)+\widetilde{g}\left(|u|^2T\, , \widetilde{\nabla}  r^{-a-b}\right)\\
& = \frac{|u|^2}{r^{a+b}}\widetilde{\Delta}r+2\left\langle\frac{uT}{r^{a+b}},\widetilde{\nabla}u\right\rangle_T-(a+b)\frac{|u|^2}{r^{a+b+1}}.
\end{aligned}
$$ Hence,  by the divergence theorem,
\[
\begin{aligned}
\frac 12\bigg |\int_{\Omega}\frac{|u|^2}{r^{a+b+1}}\left(r\widetilde{\Delta}r-(a+b)\right)dV_{\widetilde{g}}\bigg | & =
\bigg |\int_{\Omega}\left\langle\frac{uT}{r^{a+b}},\widetilde{\nabla}u\right\rangle_TdV_{\widetilde{g}}\bigg |\\
&  =
\bigg |\int_{\Omega}\left\langle\frac{uT}{r^{a}},\frac{\widetilde{\nabla}u}{r^{b}}\right\rangle_TdV_{\widetilde{g}}\bigg |,
\end{aligned}
\]
(cf. \cite[p.409]{WL}). Now applying the H\"older inequality to the right side of the above formula we obtain the desired inequalities.
\end{proof}

\noindent
{\bf Local and Global Geometric Inequalities 4.2}  {\it \, Let  $\Omega\subset M\backslash{\rm Cut}(x_0)$ be a domain in a Finsler manifold $(M, F).$ For every $u\in C_{0}^{\infty}(\Omega\backslash\{x_0\}),$
and every $a,b\in
\mathbb{R},$ the following inequality holds:
\[
\frac 1p\bigg |\int_{\Omega}\frac{|u|^p}{r^{a+b+1}}\big (r\widetilde{\Delta}r-(a+b)\big )dV_{\widetilde{g}} \bigg |\leqslant \left(\int_{\Omega}\frac{|u|^{p}}{r^{aq}}dV_{\widetilde{g}}\right)^{\frac {1}{q}}
\left(\int_{\Omega}\frac{|\widetilde{\nabla}u|^{p}_{\widetilde{g}}}{r^{bp}}dV_{\widetilde{g}}\right)^{\frac {1}{p}},
\tag{4.2.a}
\]where $p \in (1, \infty)$ and
$\frac 1p+\frac 1q=1$. In particular, if $\operatorname{Cut}(x_0)$
is empty, then}
\[
\frac 1p\bigg |\int_{M}\frac{|u|^p}{r^{a+b+1}}\big (r\widetilde{\Delta}r-(a+b)\big )dV_{\widetilde{g}} \bigg |\leqslant \left(\int_{M}\frac{|u|^{p}}{r^{aq}}dV_{\widetilde{g}}\right)^{\frac {1}{q}}
\left(\int_{M}\frac{|\widetilde{\nabla}u|^{p}_{\widetilde{g}}}{r^{bp}}dV_{\widetilde{g}}\right)^{\frac {1}{p}}.
\tag{4.2.b}
\]
\begin{proof}
First consider the case that $1<p<2$. For every $u\in
C_0^{\infty}(\Omega\backslash \left\{ x_{0}\right\})$, given $\epsilon>0$, consider $I:=p\underset{\Omega}{\displaystyle\int }\left\langle
(|u|^2+\epsilon)^\frac{p-2}{2}u\dfrac{T}{r^{a+b}},\widetilde{\nabla} u \right\rangle_{T} dV_{\widetilde{g}}\, $. Then it follows from the Guass lemma that
\[
I=\underset{\Omega}{\int }\widetilde{\rm div}\left( (|u|^2+\epsilon)^\frac{p}{2} \frac{T }{r^{a+b}}
\right) dV_{\widetilde{g}}
-\underset{\Omega}{\int
}\frac{\widetilde{\rm div}\left(T \right)}{r^{a+b}}(\left\vert
u\right\vert^2+\epsilon) ^\frac{p}{2}
dV_{\widetilde{g}}+\int_{\Omega}\frac{a+b}{r^{a+b+1}}(\left\vert
u\right\vert^2+\epsilon) ^\frac{p}{2}dV_{\widetilde{g}}.\tag{4.3}
\]
Hence by the Divergence Theorem and \eqref{31},
\begin{equation}
\begin{aligned}\label{4.4}
I &= \int_{\partial V}\left\langle
\frac{T}{r^{a+b}} (|u|^2+\epsilon)^\frac{p}{2},\xi\right\rangle_{T} dS-\int _{\Omega}\, \frac{
r\widetilde{\Delta}r-(a+b) }{r^{a+b+1}}(|u|^2+\epsilon)^\frac{p}{2}dV_{\widetilde{g}} \\
& =\epsilon^\frac{p}{2} \int_{\partial V}\left\langle
\frac{T}{r^{a+b}},\xi\right\rangle_{T} dS-\int _{\Omega}\, \frac{
r\widetilde{\Delta}r-(a+b) }{r^{a+b+1}}(|u|^2+\epsilon)^\frac{p}{2}dV_{\widetilde{g}},
\end{aligned}\tag{4.4}
\end{equation}
where $V$ is an open set with $supp\left\{ u\right\} \subset
V\subset \subset \Omega\backslash \left\{ x_{0}\right\}$, $\xi$ is the
outward unit normal vector of $\partial V$, $dS$ is the area
element induced from $dV_{\widetilde{g}}$.

\noindent
Now the triangle inequality, \eqref{4.4}, and the H\"{o}lder inequality imply that
\begin{equation}
\begin{aligned}\label{4.5}
&\left|\int_\Omega\frac{(|u|^2+\epsilon)^\frac{p}{2}}{r^{a+b+1}}(r \widetilde{\Delta} r - a - b) dV_{\widetilde{g}}\right|-\left|\epsilon^\frac{p}{2} \int_{\partial V}\left\langle
\frac{T}{r^{a+b}},\xi\right\rangle_{T} dS\right|\\
\leqslant & \left| \epsilon^\frac{p}{2} \int_{\partial V}\left\langle
\frac{T}{r^{a+b}},\xi\right\rangle_{T} dS - \int_\Omega\frac{(|u|^2+\epsilon)^\frac{p}{2}}{r^{a+b+1}}(r \widetilde{\Delta} r - a - b) dV_{\widetilde{g}}\right|\\
= & |I| \\
\leqslant& p\left( \int _{\Omega}\, \frac{\left\vert(\left\vert
u\right\vert^2+\epsilon)^\frac{p-2}{2}u
\right\vert^\frac{p}{p-1}}{r^{a\frac{p}{p-1}}}dV_{\widetilde{g}}\right)^{\frac
{p-1}{p}}\left( \int _{\Omega}\, \frac{| \widetilde{\nabla} u|_{\widetilde{g}}^{p}}{r^{bp}}dV_{\widetilde{g}}\right)^{\frac 1p}.
\end{aligned}\tag{4.5}
\end{equation}
Since $1<p<2$, we have $(|u|^2+\epsilon)^\frac{p-2}{2}<
(|u|^2)^\frac{p-2}{2}$. Thus, for every $1<p<2\, $ with $\frac
1p+\frac 1q=1$, one has via \eqref{4.5}
\begin{equation}
\begin{aligned}\label{4.6} & \left|\int_\Omega\frac{(|u|^2+\epsilon)^\frac{p}{2}}{r^{a+b+1}}(r \widetilde{\Delta} r - a - b) dV_{\widetilde{g}}\right|\\
\leqslant &
p \left( \int _{\Omega}\, \frac{\left\vert u\right\vert
^{p}}{r^{aq}}dV_{\widetilde{g}}\right)^{\frac 1q}\left( \int _{\Omega}\, \frac{|
\widetilde{\nabla} u|_{\widetilde{g}}^{p}}{r^{bp}}dV_{\widetilde{g}}\right)^{\frac 1p} + \left|\epsilon^\frac{p}{2} \int_{\partial V}\left\langle
\frac{T}{r^{a+b}},\xi\right\rangle_{T} dS\right|.
\end{aligned}\tag{4.6}
\end{equation}

Since \eqref{4.6} holds for every sufficiently small $\epsilon > 0\, ,$ we have
\begin{equation}
\begin{aligned}\label{4.7}  \left|\int_\Omega\frac{(|u|^2+\epsilon)^\frac{p}{2}}{r^{a+b+1}}(r \widetilde{\Delta} r - a - b) dV_{\widetilde{g}}\right|
\leqslant
p \left( \int _{\Omega}\, \frac{\left\vert u\right\vert
^{p}}{r^{aq}}dV_{\widetilde{g}}\right)^{\frac 1q}\left( \int _{\Omega}\, \frac{|
\widetilde{\nabla} u|_{\widetilde{g}}^{p}}{r^{bp}}dV_{\widetilde{g}}\right)^{\frac 1p}.
\end{aligned}\tag{4.7}
\end{equation}

Let $\epsilon \to 0,$ the Monotone Convergence Theorem gives the
desired result.

For the case $p\geqslant 2$, consider $I:=p\displaystyle\int_\Omega\left\langle |u|^{p-2}u\frac{T }{r^{a+b}}, \widetilde{\nabla}
u\right\rangle_T dV_{\widetilde{g}}\, $. Then it follows from the Guass lemma and \eqref{31} that
\begin{equation}
\begin{aligned}\label{4.8}I=&\int_\Omega\widetilde{\rm div}\left( \frac{|u|^p}{r^{a+b}}T \right)  dV_{\widetilde{g}}-\int_\Omega \frac{\widetilde{\rm div}\left(T\right)}{r^{a+b}}|u|^{p} dV_{\widetilde{g}} +\int_{\Omega}\frac{a+b}{r^{a+b+1}}|u|^{p}dV_{\widetilde{g}}\\
=&\int_\Omega\widetilde{\rm div}\left( \frac{|u|^p}{r^{a+b}} T
\right) dV_{\widetilde{g}}-\int_{\Omega}\frac{|u|^p}{r^{a+b+1}}\big (r\widetilde{\Delta}r-(a+b)\big )dV_{\widetilde{g}}
\end{aligned}\tag{4.8}
\end{equation}
for every $u\in C_0^\infty(\Omega\backslash \left\{ x_{0}\right\})$. Hence
by the Divergence Theorem and the H\"{o}lder inequality we easily get the desired result.
\end{proof}

\section{Proof of Theorems 1.1 and 1.2 - Generalized Caffarelli-Kohn-Nirenberg type inequalities on Finsler manifolds}

\begin{proof}[Proof of Theorem 1.1]
Case $($i$)$ $0\leqslant {\bf Ric}(\nabla r)$ and $n \leqslant a+b+1$ :

\noindent
In view of Laplacian Comparison Theorem 3.1.(2),
\begin{equation}\label{91}
\begin{aligned}
 \frac{a+b+1-n}{2}\int_{\Omega}\frac{|u|^2}{r^{a+b+1}} dV_{\widetilde{g}} & \leqslant \frac 12 \bigg |\int_{\Omega}\frac{|u|^2}{r^{a+b+1}}\left(r\widetilde{\Delta}r-(a+b)\right)dV_{\widetilde{g}} \bigg |\\
& \leqslant \left(\int_{\Omega}\frac{|u|^{p}}{r^{ap}}dV_{\widetilde{g}}\right)^{\frac {1}{p}}
\left(\int_{\Omega}\frac{|\widetilde{\nabla}u|^{q}_{\widetilde{g}}}{r^{bq}}dV_{\widetilde{g}}\right)^{\frac {1}{q}}\\
& \leqslant \left(\int_{\Omega}\frac{|u|^{p}}{r^{ap}}dV_{\rm max}\right)^{\frac {1}{p}}
\left(\int_{\Omega}\frac{|\widetilde{\nabla}u|^{q}_{\widetilde{g}}}{r^{bq}}dV_{\rm max}\right)^{\frac {1}{q}}\\
& \leqslant \mu_F^{\frac 12}\left(\int_{\Omega}\frac{|u|^{p}}{r^{ap}}dV_{\rm max}\right)^{\frac {1}{p}}
\left(\int_{\Omega}\frac{(F({\nabla}u))^{q}}{r^{bq}}dV_{\rm max}\right)^{\frac {1}{q}}.\\
\end{aligned}
\end{equation}
The second, third and last steps follow from Geometric Inequality 4.1.(a), \eqref{33}, and \eqref{32} respectively. On the other hand, \eqref{33} implies that
\begin{equation}\label{92}
\begin{aligned}
 \frac{a+b+1-n}{2}\mu_F^{-\frac {n}{2}}\int_{\Omega}\frac{|u|^2}{r^{a+b+1}} dV_{\rm {max}} & \leqslant \frac{a+b+1-n}{2}\int_{\Omega}\frac{|u|^2}{r^{a+b+1}} dV_{\rm {min}}\\
& \leqslant \frac{a+b+1-n}{2}\int_{\Omega}\frac{|u|^2}{r^{a+b+1}} dV_{\widetilde{g}}.
\end{aligned}
\end{equation}
Combining \eqref{91} and \eqref{92}, we have proved \eqref{11} when the extreme volume form $dV_{\rm {ext}}$ is $dV_{\rm {max}}$.
Similarly, we can prove \eqref{11} when  $dV_{\rm {ext}}$ is the minimum volume form $dV_{\rm {min}}$ and hence \eqref{12} holds if in addition, $\operatorname{Cut}(x_0)$
is empty or $x_0$ is a pole. This completes the proof of Case $($i$)$. Similarly, by considering the fact that a simply connected flat Finsler manifold does not have a cut point in Case $($ii$)$, and so does a simply connected Finsler manfold with nonpositive flag curvature by a comparison theorem in Case $($iii$)$, the assertions follow.
\end{proof}

\begin{proof}[Proof of Theorem 1.2]  Proceeding as in the proof of Theorem 1.1 by applying Theorem 3.1.(2), \eqref{32}, and \eqref{33} to Geometric Inequality 4.2.(a), the assertions follow.\end{proof}

\section{Proof of Theorem 1.3 }

{\rm We will give two methods:

\noindent
{\bf First Method:} This follows at once from substituting $p=q=2$ into Theorem $1.1$ or Theorem $1.2.$

\noindent
{\bf Second Method:} We first follow \cite{WL}. For every $u\in W^{1,2}_0(\Omega\backslash\{x_0\}, \widetilde{g})$ (cf. Denfinition 3.1, where $p=2$) and any $a,b,t\in\mathbb{R}$, we have, analogous to \cite[(4.1)]{WL}
$$\int_{\Omega}\left\langle\frac{\widetilde{\nabla}u}{r^b}+t\frac{u}{r^a}T,\frac{\widetilde{\nabla}u}{r^b}+t\frac{u}{r^a}T\right\rangle_TdV_{\widetilde{g}}
\geqslant0,$$namely,
\begin{equation}\label{51}\int_{\Omega}\frac{|\widetilde{\nabla}u|^2_{\widetilde{g}}}{r^{2b}}dV_{\widetilde{g}}+t^2\int_{\Omega}\frac{|u|^2}{r^{2a}}
dV_{\widetilde{g}}+2t\int_{\Omega}\left\langle\frac{uT}{r^{a+b}},\widetilde{\nabla}u\right\rangle_TdV_{\widetilde{g}}\geqslant0.
\end{equation} Observing that
$$
\begin{aligned}\widetilde{\rm div}\left(\frac{|u|^2}{r^{a+b}}T\right) & =\frac{|u|^2}{r^{a+b}}\widetilde{\Delta}r+2\widetilde{g}\left(\frac{u}{r^{a+b}}T,  \widetilde{\nabla}u\right)+\widetilde{g}\left(|u|^2T, \widetilde{\nabla} r^{-a-b}\right)\\
&=\frac{|u|^2}{r^{a+b}}\widetilde{\Delta}r+2\left\langle\frac{uT}{r^{a+b}},\widetilde{\nabla}u\right\rangle_T-(a+b)\frac{|u|^2}{r^{a+b+1}},
\end{aligned}
$$ we have by divergence theorem,
\begin{equation}\label{52}2\int_{\Omega}\left\langle\frac{uT}{r^{a+b}},\widetilde{\nabla}u\right\rangle_TdV_{\widetilde{g}}
=-\int_{\Omega}\frac{|u|^2}{r^{a+b+1}}\left(r\widetilde{\Delta}r-(a+b)\right)dV_{\widetilde{g}}.
\end{equation}Let
$$A=\int_{\Omega}\frac{|u|^2}{r^{2a}}
dV_{\widetilde{g}},\quad B=2\int_{\Omega}\left\langle\frac{uT}{r^{a+b}},\widetilde{\nabla}u\right\rangle_TdV_{\widetilde{g}},\quad
C=\int_{\Omega}\frac{|\widetilde{\nabla}u|^2_{\widetilde{g}}}{r^{2b}}dV_{\widetilde{g}},$$ then \eqref{51} takes the form
$$At^2+Bt+C\geqslant0\, ,\,  A > 0$$ for every $t\in\mathbb{R}$ which implies that $B^2-4AC\leqslant0$. Thus by \eqref{52} one has
\begin{equation}\label{53}\frac12\left|\int_{\Omega}\frac{|u|^2}{r^{a+b+1}}\left(r\widetilde{\Delta}r-(a+b)\right)dV_{\widetilde{g}}\right|\leqslant
\left(\int_{\Omega}\frac{|u|^2}{r^{2a}}dV_{\widetilde{g}}\right)^{\frac12}
\left(\int_{\Omega}\frac{|\widetilde{\nabla}u|^2_{\widetilde{g}}}{r^{2b}}dV_{\widetilde{g}}\right)^{\frac12}.\end{equation}
Suppose now that $M$ has nonnegative radial Ricci curvature at $x_0$. Then  it follows from Theorem 3.1 that $r\widetilde{\Delta}r\leqslant n-1$. Thus if in addition, $n\leqslant a+b+1$, then \eqref{53} becomes
\begin{equation}\label{54} \frac12\int_\Omega\frac{|u|^2}{r^{a+b+1}}\bigg( a+b+1-n \bigg )dV_{\widetilde{g}}\leqslant
\left(\int_\Omega\frac{|u|^2}{r^{2a}}dV_{\widetilde{g}}\right)^{\frac12}
\left(\int_\Omega\frac{|\widetilde{\nabla}u|^2_{\widetilde{g}}}{r^{2b}}dV_{\widetilde{g}}\right)^{\frac12}.\end{equation}Hence, \eqref{19} follows directly from \eqref{32}, \eqref{33} and \eqref{54}. This proves the case $($i$)$. Analogously, we can prove the cases $($ii$)$ and $($iii$)$. In particular, if $M$ has a pole $x_0$ or $\operatorname{Cut}(x_0)$ is empty in $($i$)$, or $M$ is simply connected in $($ii$)$ or $($iii$)$, then
we can choose $\Omega=M$ and \eqref{19} becomes \eqref{110}.
This completes the proof.

\noindent
{\bf Remark:} The case  that $M = \mathbb R^n$ is due to Costa $($c.f. \cite {C}$)$, The case that $M$ is a Riemannian manfiold is due to Wei-Li. $($\cite{WL}$).$

\section{Applications - Embedding Theorems for Weighted Sobolev Spaces and Differential-Integral Inequalities on Finsler manifolds}

As in the Riemannian case \cite{WL}, for giving Finsler manifold $(M,F)$ we let $L^2_{F,a}(M)$ be the completion of $C^\infty_{0}(M\backslash\{x_0\})$ with respect to the norm
$$ \|u\|_{L^2_{F,a}(M)}:=\left(\int_M\dfrac{|u|^2}{r^{2a}}dV_{\rm ext}\right)^{\frac12}, $$ $D^{1,2}_{F}(M)$ be the completion of $C^\infty_{0}(M\backslash\{x_0\})$ with respect to the  ``norm"
$$ \|u\|_{D^{1,2}_{F}(M)}:=\left(\int_M(F(\nabla u))^2dV_{\rm ext}\right)^{\frac12},$$
and $H^1_{F,a,b}(M)$ be the completion of $C^\infty_{0}(M\backslash\{x_0\})$ with respect to the  ``norm"
$$ \|u\|_{H^1_{F,a,b}(M)}:=\left(\int_M\left[\dfrac{|u|^2}{r^{2a}}+\dfrac{(F(\nabla u))^2}{r^{2b}}\right]dV_{\rm ext}\right)^{\frac12}.$$
It should be pointed out here that in general  $\|k\cdot u\|_{H^1_{F,a,b}(M)}=|k|\cdot\|u\|_{H^1_{F,a,b}(M)}$ holds only when $k\geqslant0$, thus $\|\cdot\|_{H^1_{F,a,b}(M)}$ although satisfies the triangle inequality is not a genuine norm. Neverthless, by Theorems 1.3  we clearly have

\noindent
{\bf Theorem 7.1} {\it  Let $(M,F)$ be an $n$-dimensional  complete Finsler manifold with nonpositive radial flag curvature or nonnegative radial Ricci curvature at the pole $x_0\in M$. Suppose also that $M$ has finite uniformity constant $\mu_F$.  Then the following embeddings hold

\begin{equation}
H_{F,a,b}^{1}(M) \subset L_{F, \frac {a+b+1}{2}}^{2}(M)\qquad and \qquad H_{F, b,a}^{1}(M) \subset L_{F, \frac {a+b+1}{2}}^{2}(M).\label{5.7}
\end{equation}}

As a consequence, we have  differential-integral inequalities on Finsler manifolds:

\noindent
{\bf Theorem 7.2} {\it  Let $M$ be as in Theorem 7.1. Then

\noindent
\textbf{i$)$ }For any $u\in D^{1,2}_{F}(M),$

\begin{equation}\label{5.99}
\underset{M}{\int }\frac{|u|^2}{r^{2}}dV_{\rm ext}\leqslant\left(\frac {2}{n-2}\right)^2 \cdot \mu_F^{n+1} \underset{M}{\int }|F(\nabla u)|^2 dV_{\rm ext};
\end{equation}

\noindent
\textbf{ii$)$ }For any $u\in H_{F, b+1,b}^{1}(M),$

\begin{equation}\label{5.8}
\underset{M}{\int }\frac{|u|^2}{r^{2(b+1)}}dV_{\rm ext}\leqslant\left(\frac n2-(b+1)\right)^{-2} \cdot \mu_F^{n+1} \underset{M}{\int }\frac{|F(\nabla u)|^2}{r^{2b}}dV_{\rm ext};
\end{equation}

\noindent
\textbf{iii$)$ }For any $u\in H_{F, a+1,a}^{1}(M),$

\begin{equation}\label{5.10}
\left( \underset{M}{\int }\frac{\left\vert
u\right\vert ^{2}}{r^{2(a+1)}}dV_{\rm ext}\right)^2\leqslant \left(\frac n2-(a+1)\right)^{-2} \cdot \mu_F^{n+1}\left( \underset{M}{\int }\frac{%
\left\vert u\right\vert ^{2}}{r^{2a}}dV_{\rm ext}\right) \left( \underset%
{M}{\int }\frac{\left\vert F(\nabla u)\right\vert ^{2}}{r^{2(a+1)}}dV_{\rm ext}\right);
\end{equation}

\noindent
\textbf{iv$)$ }If $u\in H_{F,-(b+1),b}^{1}(M)$ then $u$ $\in L_F^{2}(M)$
and

\begin{equation}\label{5.12}
\left( \underset{M}{\int }\left\vert u\right\vert
^{2}dV_{\rm ext}\right) ^2 \leqslant\frac{4}{n^2} \cdot \mu_F^{n+1} \left( \underset{M}{\int }r^{2(b+1)}\left\vert u\right\vert
^{2}dV_{\rm ext}\right) \left( \underset{M}{\int }\frac{\left\vert
F(\nabla u)\right\vert ^{2}}{r^{2b}}dV_{\rm ext}\right) ;
\end{equation}

\noindent
\textbf{v$)$ } If $u\in H_{F,0,1}^{1}(M),$ then $u$ $\in L_{F,1}^{2}(M)$
and

\begin{equation}\label{5.14}
\left( \underset{M}{\int }\frac{\left\vert u\right\vert ^{2}}{r^{2}}%
dV_{\rm ext}\right)^2\leqslant \frac{4}{(n-2)^2} \cdot \mu_F^{n+1} \left( \underset{M}{\int }\left\vert u\right\vert ^{2}dV_{\rm ext}\right)\left( \underset{M}{\int }\frac{\left\vert F(\nabla u)\right\vert
^{2}}{r^{2}}dV_{\rm ext}\right);
\end{equation}

\noindent
\textbf{vi$)$ }If $u\in H_{F,-1,1}^{1}(M),$ then $u$ $\in
L_{F,\frac{1}{2}}^{2}(M)$ and

\begin{equation}\label{5.16}
\left( \underset{M}{\int }\frac{\left\vert u\right\vert
^{2}}{r}dV_{\rm ext}\right) ^2 \leqslant \frac{4}{(n-1)^2} \cdot \mu_F^{n+1} \left( \underset{M}{\int }r^{2}\left\vert u\right\vert
^{2}dV_{\rm ext}\right) \left( \underset{M}{\int }\frac{\left\vert
F(\nabla u)\right\vert ^{2}}{r^{2}}dV_{\rm ext}\right) ;
\end{equation}

\noindent
\textbf{vii$)$ }If $u\in H_F^{1}(M)=$ $H_{F,0,0}^{1}(M),$ then $u$ $\in L_{F,\frac{1%
}{2}}^{2}(M)$ and

\begin{equation}\label{5.18}
 \left( \underset{M}{\int }\frac{\left\vert u\right\vert
^{2}}{r}dV_{\rm ext}\right)^2 \leqslant \frac{4}{(n-1)^2} \cdot \mu_F^{n+1} \left( \underset{M}{\int }\left\vert u\right\vert
^{2}dV_{\rm ext}\right) \left( \underset{M}{\int }\left\vert F(\nabla
u)\right\vert ^{2}dV_{\rm ext}\right).
\end{equation}}

\noindent
{\bf Remark}
The case that $M= \mathbb{R}^n$ is due to \cite{C}. The case that $M$ is a Riemannian manfiold is due to \cite{WL}. Item $(i)$ is a generalized Hardy's inequality. In the next Section we will discuss its generalizations.

\begin{proof} We make special choices in Theorems 1.3 as follows:

\noindent
\textbf{i)} \ \ \ Let $a=1, b=0;$

\noindent
\textbf{ii)} \ \ \ Let $a=b+1;$

\noindent
\textbf{iii)} \ \ Let $b=a+1;$

\noindent
\textbf{iv) }\ Let $a=-b-1;$

\noindent
\textbf{v)  }\ Let $a=0,b=1;$

\noindent
\textbf{vi)} \ \ Let $a=-1,b=1;$

\noindent
\textbf{vii) }  Let $a=0,b=0. $

\end{proof}

\section{Proof of Theorems 1.4 and 1.5 - Generalized Hardy inequalities on Finsler manifolds}
Employing the {\it double} limiting technique in [7] , we prove the following:

\noindent
{\bf Geometric Inequality 8.1} (cf. \cite[(1.3)]{WL}, \cite[(3)]{CLW2}) {\it Let $M$ be a Finsler manifold,  $u\in C_{0}^{\infty}(M)$ and $\partial B_{\delta}(x_0)$ be the $C^1$ boundary of the geodesic ball $B_{\delta}(x_0)$ centered
at $x_0$ with radius $\delta >0$. Let $V$ be an open set with smooth boundary $\partial V$ such that $V\subset\subset M$, and $u=0$ off $V$. We choose a sufficiently small $\delta>0$ so that $\partial V\cap\partial B_{\delta}(x_0)=\emptyset$. Then  for every $\epsilon>0$ and  $p>1$, we have
$$\left|-\int_{V \cap \partial B_{\delta}(x_0)}\frac{r}{r^p+\epsilon}|u|^p\langle T,\xi\rangle_T dS
+\int_{\breve{M}\backslash B_{\delta}(x_0)}\frac{(r^p+\epsilon)(1+r\widetilde{\Delta}r)-pr^p}{(r^p+\epsilon)^2}|u|^pdV_{\widetilde{g}}\right|$$
\begin{equation}\label{44} \leqslant p\left(\int_{\breve{M}\backslash B_{\delta}(x_0)}\left(\frac{|u|^{p-1}r}{r^p+\epsilon}\right)^{\frac{p}{p-1}}dV_{\widetilde{g}}\right)^{\frac{p-1}{p}}
\left(\int_{\breve{M}\backslash B_{\delta}(x_0)}|\widetilde{\nabla}u|^p_{\widetilde{g}}dV_{\widetilde{g}}\right)^{\frac{1}{p}},\end{equation}
where $\breve{M}$ is as in \eqref{3.0}, $\xi$ is the outward unit normal vector field on $\partial B_{\delta}(x_0)$ with respect to $\widetilde{g}$, and $dS$ is the volume form on $\partial B_{\delta}(x_0)$ induced from $dV_{\widetilde{g}}$.}

\begin{proof}

Observing via Gauss Lemma
$$\begin{aligned}\widetilde{\rm div}\left(\frac{rT}{r^p+\epsilon}|u|^p\right)&=\frac{\widetilde{\rm div}(rT)}{r^p+\epsilon}|u|^p+\widetilde{g}\left(rT, \widetilde{\nabla}\left(\frac{|u|^p}{r^p+\epsilon}\right)\right) \\
&=\frac{\widetilde{\rm div}(rT)}{r^p+\epsilon}|u|^p-\left\langle rT, \frac{\widetilde{\nabla} r^p}{(r^p+\epsilon)^2}|u|^p\right\rangle_T +\frac{r}{r^p+\epsilon}\langle T, \widetilde{\nabla}|u|^p\rangle_T  \\
&=\frac{\widetilde{\rm div}(rT)}{r^p+\epsilon}|u|^p-\frac{pr^p}{(r^p+\epsilon)^2}|u|^p+\frac{pru|u|^{p-2}}{r^p+\epsilon}\langle T, \widetilde{\nabla}u\rangle_T,
\end{aligned} $$ we have
\begin{equation}\label{41}\begin{aligned}&p\int_{\breve{M}\backslash B_{\delta}(x_0)}\left\langle|u|^{p-2}u\frac{rT}{r^p+\epsilon},\widetilde{\nabla}u\right\rangle_T dV_{\widetilde{g}}\\
=&\int_{\breve{M}\backslash B_{\delta}(x_0)}\widetilde{\rm div}\left(\frac{rT}{r^p+\epsilon}|u|^p\right)dV_{\widetilde{g}}
-\int_{\breve{M}\backslash B_{\delta}(x_0)}\frac{\widetilde{\rm div}(rT)}{r^p+\epsilon}|u|^pdV_{\widetilde{g}}
+\int_{\breve{M}\backslash B_{\delta}(x_0)} \frac{pr^p}{(r^p+\epsilon)^2}|u|^pdV_{\widetilde{g}}.\end{aligned}\end{equation} By the divergence theorem it follows that for sufficiently small $\delta > 0$,
\begin{equation}\begin{aligned}\label{42}\int_{\breve{M}\backslash B_{\delta}(x_0)}\widetilde{\rm div}\left(\frac{rT}{r^p+\epsilon}|u|^p\right)dV_{\widetilde{g}}
&=\int_{V\backslash B_{\delta}(x_0)}\widetilde{\rm div}\left(\frac{rT}{r^p+\epsilon}|u|^p\right)dV_{\widetilde{g}}\\
&=\int_{V \cap \partial B_{\delta}(x_0)}\frac{r}{r^p+\epsilon}|u|^p\langle T,\xi\rangle_T dS.\end{aligned}\end{equation}  On the other hand,
\begin{equation}\label{43} \widetilde{\rm div}(rT)=r\widetilde{\rm div}(T)+\widetilde{g}(T, \widetilde{\nabla} r)=1+r\widetilde{\Delta}r.\end{equation}
Substituting \eqref{42} and \eqref{43} into \eqref{41} and using the H\"{o}lder inequality one has
the desired \eqref{44}.
\end{proof}

\begin{proof}[Proof of Theorem 1.4]
We first assume $u \in C_0^{\infty}(M)\, .$ Since $1<p<n$ and $M$ has nonpositive radial flag curvature at $x_0$, by Theorem 3.1 we have $r\widetilde{\Delta}r+1\geqslant n>p$, and thus
\begin{equation}\label{8.5} \begin{aligned} \int_{\breve{M}\backslash B_{\delta}(x_0)}\frac{(r^p+\epsilon)(1+r\widetilde{\Delta}r)-pr^p}{(r^p+\epsilon)^2}|u|^pdV_{\widetilde{g}}
& \geqslant \int_{\breve{M}\backslash B_{\delta}(x_0)}\frac{(n-p)r^p+(n-p)\epsilon}{(r^p+\epsilon)^2}|u|^pdV_{\widetilde{g}}\\
& \geqslant(n-p)\int_{\breve{M}\backslash B_{\delta}(x_0)}\frac{(r^p+\epsilon)^{\frac{1}{p-1}}}{(r^p+\epsilon)^{\frac{p}{p-1}}}|u|^pdV_{\widetilde{g}}\\
& \geqslant(n-p)\int_{\breve{M}\backslash B_{\delta}(x_0)}\frac{(r^p)^{\frac{1}{p-1}}}{(r^p+\epsilon)^{\frac{p}{p-1}}}|u|^pdV_{\widetilde{g}}.
\end{aligned}\end{equation}
Substituting \eqref{8.5} into \eqref{44},  we have for sufficiently small $\delta>0$, 
$$-\int_{V\cap\partial B_{\delta}(x_0)}\frac{r}{r^p+\epsilon}|u|^p\langle T,\xi\rangle_T dS
+(n-p)\int_{\breve{M}\backslash B_{\delta}(x_0)}\frac{(r^p)^{\frac{1}{p-1}}}{(r^p+\epsilon)^{\frac{p}{p-1}}}|u|^pdV_{\widetilde{g}}$$
\begin{equation}\label{45} \leqslant p\left(\int_{\breve{M}\backslash B_{\delta}(x_0)}\frac{(r^p)^{\frac{1}{p-1}}}{(r^p+\epsilon)^{\frac{p}{p-1}}}|u|^pdV_{\widetilde{g}}\right)^{\frac{p-1}{p}}
\left(\int_{\breve{M}\backslash B_{\delta}(x_0)}|\widetilde{\nabla}u|^p_{\widetilde{g}}dV_{\widetilde{g}}\right)^{\frac{1}{p}}.\end{equation}
For sufficiently small $\delta>0$, one has
\begin{equation}\label{8.7}\int_{\partial B_{\delta}(x_0)}\frac{r}{r^p+\epsilon}|u|^p\langle T,\xi\rangle_T dS=0\quad \text{if}\quad x_0\notin V\end{equation} and
\begin{equation}\label{8.8}\left|\int_{\partial B_{\delta}(x_0)}\frac{r}{r^p+\epsilon}|u|^p\langle T,\xi\rangle_T dS\right|
\to 0\quad \text{as}\quad \delta\to 0\quad \text{if}\quad x_0\in V,\end{equation}
Indeed,  $\frac{r}{r^p+\epsilon}$ is a continuous, nondecreasing function for $r \in [0, \delta _0]$, where $\delta _0 = (\frac{\epsilon}{p-1})^{\frac 1p}$ and $u$ is bounded in $M$. Hence,  for $\delta < \delta _0\, ,$
\begin{equation}\label{8.9}\left|\int_{\partial B_{\delta}(x_0)}\frac{r}{r^p+\epsilon}|u|^p\langle T,\xi\rangle_T dS\right|
\leqslant\frac{\delta}{\delta ^p+\epsilon} \int _{\partial B_{\delta}(x_0)} \max_M |u|^p dS.
\end{equation} This implies \eqref{8.8}.
It follows from \eqref{45}, via \eqref{8.7} or \eqref{8.8} 
that for every $\epsilon>0$,
\begin{equation}\label{46}\begin{aligned} (n-p)\left(\int_{\breve{M}\backslash B_{\delta}(x_0)}\frac{(r^p)^{\frac{1}{p-1}}}{(r^p+\epsilon)^{\frac{p}{p-1}}}|u|^pdV_{\widetilde{g}}\right)^{\frac{1}{p}}&\leqslant p
\left(\int_{\breve{M}\backslash B_{\delta}(x_0)}|\widetilde{\nabla}u|^p_{\widetilde{g}}dV_{\widetilde{g}}\right)^{\frac{1}{p}}\\
& \leqslant p \left(\int_{\breve{M}}|\widetilde{\nabla}u|^p_{\widetilde{g}}dV_{\widetilde{g}}\right)^{\frac{1}{p}}.
\end{aligned}\end{equation}
\noindent
Monotone convergence theorem and \eqref{33} imply that as $\epsilon \to 0$, for every $u \in C_0^{\infty}(M)$,
\begin{equation}\label{8.11} \bigg (\frac {n-p}{p}\bigg ) \left(\int_{\breve{M}}\bigg | {\frac{u}{r}}\bigg |^p dV_{\widetilde{g}}\right)^{\frac 1p}\leqslant
 \left(\int_{\breve{M}}|\widetilde{\nabla}u|^p_{\widetilde{g}}dV_{\widetilde{g}}\right)^{\frac {1}{p}},\end{equation}
\begin{equation}\label{8.12}\quad \text{and}\quad \frac {u}{r} \in L^p(M, \widetilde{g}).
\end{equation}

Now we {\it extend \eqref{8.11} from $u \in C^{\infty}_{0} (M)$ to $u \in W_{F,0}^{1,p}(M)\, .$} Let $\{u_i\}$ be a sequence of functions in  $ C_0^{\infty}(M)$ tending to $u \in W_{F,0}^{1,p}(M)\, $ in $
\|\, \cdot \,\|_{W_{F,0}^{1,p}(M)}$ as in \eqref{1.11}.
Applying the inequality \eqref{8.11} to difference
$u_{i_m} - u_{i_n}\, ,$ we have via \eqref{32} and \eqref{33}

\begin{equation}\begin{aligned}
 \left(\int_{\breve{M}}\bigg | {\frac{|u_{i_m} - u_{i_n}}{r}}\bigg |^p dV_{\widetilde{g}}\right)^{\frac 1p}&\leqslant \bigg (\frac {p}{n-p}\bigg )
 \left(\int_{\breve{M}}|\widetilde{\nabla}(u_{i_m} - u_{i_n})|^p_{\widetilde{g}}dV_{\widetilde{g}}\right)^{\frac {1}{p}}\\
& \leqslant \bigg (\frac {p}{n-p}\bigg ) \mu _F^{\frac{1}{2}}
 \left(\int_{M}\bigg (F\big (\nabla(u_{i_m} - u_{i_n})\big )\bigg )^p dV_{\text{max}}\right)^{\frac {1}{p}}.\end{aligned}
 \label{8.13}
\end{equation}

\noindent
Hence $\{\frac {u_i}{r}\}$ is a Cauchy sequence  in $L^p(M, \widetilde{g})\, .$ By Theorem 3.3, there exists a limiting function $f(x) \in L^p(M,
\widetilde{g})\, $ satisfying,
 via \eqref{32} and \eqref{33},

\begin{equation}
\begin{aligned}
\int _{\breve{M}} |f(x)|^p\, dV_{\widetilde{g}} = \lim_{i\to \infty} \int _{\breve{M}} \frac {|u_i(x)|^p}{r^p}\, dV_{\widetilde{g}}
& \leqslant \bigg (\frac {p}{n-p}\bigg )^p \lim_{i\to \infty} \int _{\breve{M}} |\widetilde{\nabla} u_{i} |^p_{\widetilde{g}} \, dV_{\widetilde{g}}\\
& \leqslant \bigg (\frac {p}{n-p}\bigg )^p \mu _F^{\frac{p}{2}}\lim_{i\to \infty} \int _{M} (F(\nabla u_{i}))^p \, dV_{\text{max}}\\
& \leqslant \bigg (\frac {p}{n-p}\bigg )^p \mu _F^{\frac{p}{2}} \int _{M} (F(\nabla u))^p \, dV_{\text{max}}\\
& \leqslant \bigg (\frac {p}{n-p}\bigg )^p \mu _F^{\frac{n+p}{2}} \int _{M} (F(\nabla u))^p \, dV_{\text{min}}.
\label{8.14}
\end{aligned}
\end{equation}
On the other hand, since $\frac {1}{r^p}$ is bounded in $M\backslash B_{\epsilon}(x_0)\, ,$ where $ B_{\epsilon}(x_0)$ is the open geodesic ball of radius $\epsilon > 0\, ,$ centered at $x_0$, and the pointwise convergence in Theorem 3.4.(2), we have for every   $\epsilon > 0\, ,$

\begin{equation}
\begin{aligned}
\int _{\breve{M}\backslash B_{\epsilon}(x_0)} |f(x)|^p\, dV_{\widetilde{g}} =  \lim_{i\to \infty} \int _{\breve{M}\backslash B_{\epsilon}(x_0)} \frac {|u_i(x)|^p }{r^p}\, dV_{\widetilde{g}} &= \int _{\breve{M}\backslash B_{\epsilon}(x_0)} \frac {|u|^p }{r^p}\, dV_{\widetilde{g}}\\
& = \int _{\breve{M}} \chi _{\breve{M}\backslash B_{\epsilon}(x_0)} \frac {|u|^p }{r^p}\, dV_{\widetilde{g}},
\label{8.15}
\end{aligned}
\end{equation}
where $\chi _{\breve{M}\backslash B_{\epsilon}(x_0)}$ is the characteritic function on $ \breve{M}\backslash B_{\epsilon}(x_0)$.
As $\epsilon \to 0\, ,$ monotone convergence theorem and \eqref{33} imply that
\begin{equation}
\begin{aligned}
\int _{\breve{M}} |f(x)|^p\, dV_{\widetilde{g}} =  \lim_{i\to \infty} \int _{\breve{M}} \frac {|u_i|^p }{r^p}\, dV_{\widetilde{g}} & = \int _{\breve{M}} \frac {|u|^p }{r^p}\, dV_{\widetilde{g}}\\
&  \geqslant \int _{\breve{M}} \frac {|u|^p }{r^p}\, dV_{\text{min}}\\
&  \geqslant \mu _F^{-\frac {n}{2}}\int _{\breve{M}} \frac {|u|^p }{r^p}\, dV_{\text{max}}.
\end{aligned}
\label{8.16}\end{equation}
Substituting \eqref{8.16} into \eqref{8.14} we have  for every $u \in W_{F,0}^{1,p}(M)$,
$$\int _{\breve{M}} \frac {|u|^p }{r^p}\, dV_{\rm ext} \leqslant \bigg (\frac {p}{n-p}\bigg )^p
\mu _F^{\frac{n+p}{2}} \int _{M} (F(\nabla u))^p \, dV_{\rm ext},  $$ which implies that
$$ \int _{M} \chi _{M\backslash B_{\epsilon}(x_0)}\frac {|u|^p }{r^p}\, dV_{\rm ext}=
\int_{M\backslash B_{\epsilon}(x_0)} \frac {|u|^p }{r^p}\, dV_{\rm ext}
=\int_{\breve{M}\backslash B_{\epsilon}(x_0)} \frac {|u|^p }{r^p}\, dV_{\rm ext}  $$
$$ \leqslant \bigg (\frac {p}{n-p}\bigg )^p
\mu _F^{\frac{n+p}{2}} \int _{M} (F(\nabla u))^p \, dV_{\rm ext}, \;\forall \epsilon>0.  $$
As $\epsilon \to 0\, ,$ again by monotone convergence theorem we obtain  the desired $\frac {u}{r} \in L_{F}^p(M)$ and inequality \eqref{111}.

\end{proof}

The proof of Theorem 1.4 yields

\noindent
{\bf Corollary 8.1} {\it Let $(M,F)$ be an $n$-dimensional  complete Finsler manifold with nonpositive radial flag curvature at the pole $x_0\in M$ and with finite uniformity constant $\mu_F$. Then for any $u\in W_{0}^{1,p}(M,\tilde{g})$ and $1<p<n$, the following Hardy type inequality holds:}
\[ \bigg (\frac {n-p}{p}\bigg ) \left(\int_{\breve{M}}\bigg | {\frac{u}{r}}\bigg |^p dV_{\widetilde{g}}\right)^{\frac 1p}\leqslant
 \left(\int_{\breve{M}}|\widetilde{\nabla}u|^p_{\widetilde{g}}dV_{\widetilde{g}}\right)^{\frac {1}{p}}.\]
{\it Furthermore, $\frac {u}{r} \in L^p(M, \tilde{g})\, .$}\smallskip

Either Theorem 1.4 or Corollary 8.1 recaptures the following, when $M$ is a Riemannian manifold.\smallskip

\noindent
{\bf Corollary 8.2} (\cite[Theorem 1, Corollary 1.2]{WL}) {\it Let $M$ be an $n$-dimensional  complete Riemannian manifold of nonpositive radial curvature with the volume element $dv$. Then for any $u\in W_{0}^{1,p}(M)$ and $1 < p < n$, the following Hardy type inequality holds:
\begin{equation}\label{8.21}
\left(\frac{n-p}{p}\right)^p  \int_M\frac{|u|^p}{r^p} dv\leqslant \int_M |\nabla u|^p \, dv.
\end{equation}
Furthermore, $\frac {u}{r} \in L^p(M)\,  .$}

\begin{proof}[Proof of Theorem 1.5.] We first assume $u\in C_0^{\infty} (M)\, .$ When $p>n$ and $M$ has nonnegative radial Ricci curvature at $x_0$, by Theorem 3.1 we have $r\widetilde{\Delta}r+1\leqslant n<p$. In view of this inequality
and  the triangle inequality, \eqref{44} implies
\begin{equation}\label{8.17} \begin{aligned} \int_{\breve{M}\backslash B_{\delta}(x_0)}\frac{(p-n)r^p-n\epsilon}{(r^p+\epsilon)^2}|u|^pdV_{\widetilde{g}}
& \leqslant p\left(\int_{\breve{M}\backslash B_{\delta}(x_0)}\left(\frac{|u|^{p-1}r}{r^p+\epsilon}\right)^{\frac{p}{p-1}}dV_{\widetilde{g}}\right)^{\frac{p-1}{p}}
\left(\int_{\breve{M}\backslash B_{\delta}(x_0)}|\widetilde{\nabla}u|^p_{\widetilde{g}}dV_{\widetilde{g}}\right)^{\frac{1}{p}} \\ & \quad+ \left|\int_{\partial B_{\delta}(x_0)}\frac{r}{r^p+\epsilon}|u|^p\langle T,\xi\rangle_T dS\right|.
\end{aligned}\end{equation}
Applying \eqref{8.7} or \eqref{8.8}, and letting $\delta \to 0$ in \eqref{8.17}, one has
\begin{equation}\label{8.18}  \int_{\breve{M}}\frac{(p-n)r^p-n\epsilon}{(r^p+\epsilon)^2}|u|^pdV_{\widetilde{g}}
\leqslant p\left(\int_{\breve{M}}\left(\frac{|u|^{p-1}r}{r^p+\epsilon}\right)^{\frac{p}{p-1}}dV_{\widetilde{g}}\right)^{\frac{p-1}{p}}
\left(\int_{\breve{M}}|\widetilde{\nabla}u|^p_{\widetilde{g}}dV_{\widetilde{g}}\right)^{\frac{1}{p}}.\end{equation}
We observe that the integrands in the left, and in the first factor in the right of \eqref{8.18} are monotone and  uniformly bounded above by a positive constant multiple of $\big |\frac {u}{r}\big |^p$ on $M\, .$ Since $\frac {u}{r} \in L_F^p(M)\, , \frac {u}{r} \in L^p(M, \widetilde{g})$. By the dominated convergent theorem, as $\epsilon\to 0\, ,$
\begin{equation}\label{47}  (p-n) \left(\int_{\breve{M}}\left|\frac {u}{r}\right|^p dV_{\widetilde{g}}\right)
\leqslant p\left(\int_{\breve{M}}\left|\frac {u}{r}\right|^p dV_{\widetilde{g}}\right)^{\frac{p-1}{p}}
\left(\int_{\breve{M}}|\widetilde{\nabla}u|^p_{\widetilde{g}}dV_{\widetilde{g}}\right)^{\frac{1}{p}}.\end{equation}
Simplifying and raising to the $p$-th power,
\begin{equation}\label{48}  \left(\frac{p-n}{p}\right)^p \left(\int_{\breve{M}}\left|\frac {u}{r}\right|^p dV_{\widetilde{g}}\right)
\leqslant
\left(\int_{\breve{M}}|\widetilde{\nabla}u|^p_{\widetilde{g}}dV_{\widetilde{g}}\right).\end{equation}
Now analogously we {\it extend \eqref{48} from $u \in C^{\infty}_{0}(M)$ to $u \in W_{F,0}^{1,p}(M)\, .$} Let $\{u_i\}$ be a sequence of functions in  $ C_0^{\infty}(M)$ tending to $u \in W_{F,0}^{1,p}(M)\, .$  Applying the inequality \eqref{48} to difference
$u_{i_m} - u_{i_n}\, ,$ employing \eqref{32} and \eqref{33} and proceeding as in the proof of Theorem 1.4,   we obtain the desired inequality \eqref{112} for every $u \in W_{F,0}^{1,p}(M)\, .$
\end{proof}

Similarly, the proof of Theorem 1.4 yields

\noindent
{\bf Corollary 8.3} {\it Let $(M,F)$ be an $n$-dimensional  complete Finsler manifold with nonnegative radial Ricci curvature at the pole $x_0\in M$, and with finite uniformity constant $\mu_F$. Then for any $u\in W_{0}^{1,p}(M, \tilde{g}),$ $\frac {u}{r} \in L^p(M,\tilde{g})$  and $p>n$, the following Hardy type inequality holds:}
\[  \left(\frac{p-n}{p}\right)^p \left(\int_{\breve{M}}\left|\frac {u}{r}\right|^p dV_{\widetilde{g}}\right)
\leqslant
\left(\int_{\breve{M}}|\widetilde{\nabla}u|^p_{\widetilde{g}}dV_{\widetilde{g}}\right).\]
\smallskip

Either Theorem 1.5 or Corollary 8.3 recaptures the following, when $M$ is a Riemannian manifold.\smallskip

\noindent

\noindent
{\bf Corollary 8.4} (\cite[Theorem 5]{CLW2}) {\it Let $M$ be an $n$-dimensional  Riemannian manifold with a pole, nonnegative radial Ricci curvature and the volume element $dv$. Then for any $u\in W_{0}^{1,p}(M)\, ,$ $\frac {u}{r} \in L^p(M)$ and $p>n$, the following Hardy type inequality holds:}
\begin{equation}\label{8.22}
\left(\frac{p-n}{p}\right)^p \int_M\frac{|u|^p}{r^p} dv\leqslant \int_M |\nabla u|^p \, dv.
\end{equation}

\section*{Acknowledgments}  

\noindent
$^*$ The first author is supported in part by NSF (DMS-1447008), and the OU Arts and Sciences Travel Assistance Program;\\
$^{**}$ The corresponding author is  supported in part by the National Science Foundation of China (No. 12001259) and the National Science Foundation of Fujian province of China (No. 2020J01131142).\\
The authors wish to thank the referees for their helpful comments and suggestions.

\end{document}